\documentclass[letterpaper,12pt]{amsart}

\usepackage[dvipsnames]{xcolor}
\usepackage{latexsym,array,delarray,amsthm,amssymb,epsfig,setspace,tikz,amsmath,enumerate,mathrsfs,graphicx,mathtools}
\usepackage{ wasysym }
\usepackage[shortlabels]{enumitem}
\usepackage[ruled,vlined]{algorithm2e}
\usepackage{blkarray}
\usepackage{url}
\usepackage{microtype}
\usepackage[colorlinks=true,hyperindex, linkcolor=magenta, pagebackref=false, citecolor=cyan]{hyperref}
\usepackage[alphabetic,lite,backrefs]{amsrefs}
\usepackage{subcaption}
\usepackage{cleveref}
\usepackage[backgroundcolor=white,bordercolor=black]{todonotes}
\usepackage{comment}

\usepackage{tikz-cd} 
\usetikzlibrary{arrows,automata}

\theoremstyle{plain}
\newtheorem{thm}{Theorem}[section]
\newtheorem{lm}[thm]{Lemma}
\newtheorem{prop}[thm]{Proposition}
\newtheorem{cor}[thm]{Corollary}

\theoremstyle{definition}
\newtheorem{defn}[thm]{Definition}
\newtheorem{ex}[thm]{Example}

\newcommand{\Z}{\mathbb{Z}}
\newcommand{\N}{\mathbb{N}}

\newcommand{\K}{\mathbb{K}}

\newcommand{\ba}{\mathbf{a}}
\newcommand{\bb}{\mathbf{b}}

\newcommand{\bk}{\mathbf{k}}
\newcommand{\bl}{\mathbf{l}}

\newcommand{\bu}{\mathbf{u}}
\newcommand{\bv}{\mathbf{v}}

\newcommand{\bx}{\mathbf{x}}

\newcommand{\bm}{\mathbf{m}}

\newcommand{\bs}{\mathbf{s}}

\newcommand{\btau}{\boldsymbol{\tau}}

\newcommand{\cA}{\mathcal{A}}

\newcommand{\G}{\mathcal{G}}
\newcommand{\cL}{\mathcal{L}}

\newcommand{\cN}{\mathcal{N}}

\newcommand{\sA}{\mathscr{A}}
\newcommand{\sB}{\mathscr{B}}
\newcommand{\sC}{\mathscr{C}}
\DeclareMathOperator{\equivH}{{equivH}}
\newcommand{\sh}{\mathrm{sh}}

\DeclareMathOperator{\Mon}{Mon}

\DeclareMathOperator{\FI}{FI}
\DeclareMathOperator{\OI}{OI}
\DeclareMathOperator{\Sym}{Sym}
\DeclareMathOperator{\Inc}{Inc}
\DeclareMathOperator{\GL}{GL}

\begin{document}

\title[Segre Products and Regular Languages]{Shift Invariant Algebras, Segre Products and Regular Languages }

\author{Aida Maraj}
\address{Department of Mathematics\\
University of Michigan \\
1855 East Hall, Ann Arbor, MI 48109  \\
 USA}
\email{maraja@umich.edu}

\author{Uwe Nagel}
\address{Department of Mathematics\\
University of Kentucky\\
715 Patterson Office Tower\\
Lexington, KY 40506-0027 \\ USA}
\email{uwe.nagel@uky.edu}

%\author{Aida Maraj and Uwe Nagel}

\thanks{The second author was partially supported by Simons Foundation grant \#636513. }

\maketitle

\begin{abstract}
Motivated by results on the rationality of equivariant Hilbert series of some hierarchical models in algebraic statistics  we introduce the Segre product of formal languages and apply it to establish rationality of equivariant Hilbert series in new cases. To this end we show that the Segre product of two regular languages is again regular. We also prove that every filtration of algebras given as a tensor product of families of algebras with  rational equivariant Hilbert series has a rational equivariant Hilbert series. The term equivariant is used broadly to include the action of the monoid of nonnegative integers by shifting variables. Furthermore, we exhibit a filtration of  shift invariant monomial algebras that has a rational equivariant Hilbert series, but whose presentation ideals do not stabilize.  
\end{abstract}

%\tableofcontents

%%%%%%%%%%%%%%%%%%%%%%%%%%%%

%\input{Paper/Chapter_Introduction}
\section{Introduction}

In \cite{HS} (see also \cite{C}), Hillar and Sullivant initiated a systematic study of filtrations of ideals 
$(I_n)_{n \in \N}$ of ideals $I_n \subset \K[X_{[c]\times [n]}] = \K[x_{i, j} \; \mid \; 1 \le i \le c, 1 \le j \le n]$ with $c \in \N$ that 
 are rather symmetric. They showed that the ideals $I_n$ stabilize, that is, their colimit $I = {\displaystyle \lim_{\longrightarrow}}\, I_n$ in $\K[X_{[c]\times \N}] = \K[x_{i, j} \; \mid \; 1 \le i \le c, j \in \N]$ is finitely generated up to symmetry. The equivariant Hilbert series of such a filtration was introduced in \cite{NR} and shown to be rational. The stabilization result was extended to modules in \cite{NR2} by establishing that the filtration $(\K[X_{[c]\times [n]}])_{n \in \N}$ is a noetherian $\OI$- and $\FI$-algebra with the property that finitely generated modules over it are noetherian. The rationality result of equivariant Hilbert series was refined and extended to modules in \cite{N}. 
Here, $\OI$ and $\FI$ denote suitable combinatorial categories. 
%the categories whose objects are finite (ordered) sets and whose morphisms are strictly increasing and bijective maps, respectively. 

In \cite{NR2}, it was also shown that polynomial $\OI$- and $\FI$-algebras other than $(\K[X_{[c]\times [n]}])_{n \in \N}$ are not noetherian. In particular, Segre products of algebras are typically not quotients of a noetherian $\OI$- or $\FI$-algebra. Thus, there is no general result for guaranteeing  rationality of an equivariant Hilbert series of Segre products. Our main motivation is to address this problem by introducing a new method. The starting point is a rationality result on equivariant Hilbert series of some hierarchical models in algebraic statistics in \cite{MN}.  
%Nevertheless,  we established the rationality of equivariant Hilbert series of some hierarchical models in algebraic statistics by using formal languages .  
In this paper we generalize this approach considerably. The main novelty is the introduction of the Segre product of formal languages  (see \Cref{def:Segre language}). It can be used to enumerate monomials in a Segre product of monomial algebras. 
As an application we establish rationality of equivariant Hilbert series in new cases. Regular languages have been used previously to prove rationality results (see, for example, \cite{SS-14,SS-16, KLS}). 

Furthermore, we consider tensor products of algebras and show that every 2-filtration of algebras given as tensor product of 1-filtrations of algebras with a rational equivariant Hilbert series witnessed by a regular language also has a rational equivariant Hilbert series (see \Cref{prop:relate tensor products}). 

In this paper we focus on filtrations of monomial subalgebras and ideals whose colimits are invariant under shifting variables $x_{i, j}$ to $x_{i, j+1}$. These algebras and their toric presentation ideals have smaller automorphism groups than, for example, $\FI$-ideals. The latter correspond to filtrations of ideals that are invariant under the action of a symmetric group. 
We exhibit a filtration of shift invariant monomial algebras that has a rational equivariant Hilbert series, but whose presentation ideals do not stabilize.  In particular, their colimit is not finitely generated up to shifting. 

We now describe the organization of this article. In the following section we review some basic concepts. 
 Filtrations determined by Segre products of  filtrations of algebras are studied in \Cref{sec:Segre Products of Languages and Algebras}. To this end we introduce the Segre product of two formal languages (see \Cref{def:Segre language}) and show in \Cref{prop:relate Segre products} that it represents the monomials in a Segre product. We also establish 
 that the Segre product of two regular languages is again regular (see \Cref{thm:regularity of Segre language}).
 It follows that the Segre product of filtrations whose factors are represented by regular languages has a rational equivariant Hilbert series (see \Cref{prop:relate Segre products}).  The analogous questions for 2-filtrations determined by tensor products are discussed in \Cref{sec:tensor product}, where we establish the mentioned rationality result for their equivariant Hilbert series. 
  In \Cref{sec:1}, we consider an infinite family of filtrations of monomial algebras with shift invariant colimits. We use formal languages to compute the equivariant Hilbert series of any such filtration.  Thus, Segre products of these filtrations have a rational equivariant Hilbert series as well (see \Cref{thm:Segre product with rat hilb}). 
  %This is a consequence of our results on Segre products of languages. 
We conclude this paper by discussing a filtration  of monomial algebras which seems at first glance very similar to the filtrations considered in \Cref{sec:1}. Yet, it exhibits a new phenomenon. Their toric presentation ideals do not stabilize by \Cref{cor:ideal not fg}, but their equivariant Hilbert series is rational (see \Cref{cor:Hilb of non-fg ideal}). 

%%%%%%%%%%%%%%%%%%%%%%%%%%%%%%%%%%%%%%%%%%%%%%%%

\section{Basic Concepts} 
  \label{sec:prep} 
  
We review some concepts and techniques we use in subsequent sections. For more details, we refer to the textbooks \cite{BH, HHO, HU}.  

A \emph{standard graded} algebra  over a field $K$ is a graded $\K$-algebra $A = \oplus_{j \ge 0} [A]_j$ with $[A]_0 = \K$ that is generated in degree one. All polynomial rings in this note are standard graded. We set $\K[X_{[c]\times [n]}] = \K[x_{i, j} \; \mid \;  i \in [c], j \in [n]]$ with $c \in \N$, where $[n] = \{1,2,\ldots,n\}$ and $[0] = \emptyset$.  If $A$ is noetherian its \emph{Hilberts series} is a formal power series in one variable $t$:
\[
H_{A}(t)\coloneqq \sum_{j\geq 0}\dim_K[A]_{j}t^j.
\]
It is a rational function in $\mathbb{Q}(t)$. 

Often we consider a family  $\sA=(A_n)_{n\in \N}$  of noetherian standard graded algebras $A_n$. Typically, some compatibility conditions between the algebras $A_n$ are assumed. $\OI$- and $\FI$-algebras as defined in \cite{NR2} may be viewed this way.  As another example, we will consider filtrations (see \Cref{def:represented by language}).  
Following \cite{NR}, we define the \emph{equivariant Hilbert series}  of a family  $\sA=(A_n)_{n\in \N}$ as a formal power series  in two variables $s$ and $t$ as
\[
\equivH_{\sA} (t,s)\coloneqq \sum_{n\geq 1} H_{A_n}(t)s^{n} =  \sum_{n\geq 1} \sum_{j \ge 0} \dim_{\K} [A_n]_j t^j s^n. 
\]
We are interested in establishing instances where this series is a rational function. 

Using the standard embedding of $\K[X_{[c]\times [m]}]$ into $\K[X_{[c]\times [n]}]$ if $m \le n$, the colimit
$A = {\displaystyle \lim_{\longrightarrow}}\, \K[X_{[c]\times [n]}] \cong  \K[X_{[c]\times \N}] = \K[x_{i, j} \; \mid \; i \in [c], j \in \N]$ is 
 a polynomial ring in infinitely many variables. It is invariant under the action of the monoid $\N_0$ by \emph{shifting} the second index of a variable, that is, $\sh_k (x_{i, j}) = x_{i, j+k}$. We will consider shift invariant ideals and quotient algebras, complementing previous investigations of structures that are invariant under the action of groups such as $\Sym_\infty$ or $\GL_{\infty}$ or the monoid $\Inc$ of increasing functions  (see, e.g., \cite{DEKL, HS, SS-12b}). 

\begin{ex}
   \label{exa:review}
Consider a family $\sA=(A_n)_{n\in \N}$  with $A_n = \K[x_i^2,x_ix_{i+1}\:\mid\: i \in  [n] ]$. Its colimit is 
$A = \K[x_i^2,x_ix_{i+1}  \:\mid\: i \in  \N ]$. It is invariant under shifting, but it is not invariant under the action of $\Sym_\infty$ or $\Inc$. Moreover, the colimit $A=\lim\limits_{n\rightarrow \infty} A_n$ is finitely generated as an algebra by $x_1^2,x_1x_2$ up to shifting indices.

Note that each algebra $A_n$ is a quotient of $\K[X_{[2]\times [n]}]$ with  \emph{presentation ideal}  $I_n = \ker \varphi_n$, where $\varphi_n$ is the algebra homomorphism 
\[
\varphi_n \colon  \K[X_{[2]\times [n]}] \longrightarrow \K[X_{\N}] = \K[x_i \; \mid \; i \in \N], 
\]
defined by $x_{1,i} \mapsto  x_i^2, \quad x_{2,i} \mapsto x_i x_{i+1}$. Their colimit $I = {\displaystyle \lim_{\longrightarrow}}\, I_n$ is shift invariant. In fact, one can show that it is generated by $x_{1,i} x_{1,{i+1}}-x_{2,{i}}^2$ with $i \in \N$, that is, $I$ is generated by $x_{1,1}x_{1,2}-x_{2,1}^2$  up to shifting. 

Assigning each generating monomial of $A_n$ degree one, $A_n$ becomes a standard graded $\K$-algebra, and there is a graded algebra isomorphism $A_n \cong \K[X_{[2]\times [n]}]/I_n$. 
\end{ex}
  
Determining Hilbert series can be challenging. Formal languages and finite automata can help in some cases.  A \emph{formal language} $\cL$ is any subset of a free monoid $\Sigma^*$, where $\Sigma$ be a finite set. We refer to the elements of $\cL$ as \emph{words} in the alphabet $\Sigma$. The class of \emph{regular languages} on $\Sigma$ is the smallest class of languages that contains the languages having a letter of $\Sigma$ or the empty word as their only  word and that is closed under taking unions, concatenation and passing from a language $\cN$ to its Kleene star $\cN^*$ (see \cite[Section 4.2]{HU} for details and unexplained terminology). 

\begin{ex} 
    \label{exa:easy language}
Consider the language $\cL$ on the alphabet $\Sigma=\{ \alpha_{0},\alpha_{1},\dots,\alpha_{c}, \tau\}$ in  $c+2$ letters defined as 
\begin{align*}
\cL & = \{ \tau^{k_1}\alpha_{i_1}\tau^{k_2}\alpha_{i_2}\dots \tau^{k_d}\alpha_{i_d}\tau^{k_{d+1}} \; \mid \; d \in \N_0,   0 \le i_1,\dots,i_d \le c, k_1 \ge 0, \\
& \hspace*{8cm} k_{j+1 }\ge i_{j} \text{ if }  1\le j \le d \}. 
\end{align*}
The language $\cL$ is regular. Indeed, the words in $\cL$ are exactly the words in the alphabet $\Sigma'=\{ \beta_{0},\dots,\beta_{c}, \tau\}$, i.e., $\cL = (\Sigma')^*$ with $\beta_j = \alpha_j \tau^j$ for $j=0,\ldots,c$. 
\end{ex}
 
Equivalently, a regular language is a  language that can be realized by a \emph{finite automaton}. 
A \emph{finite automaton} is a labeled directed graph whose vertices represent states and edges represent transitions between 
states. It has an \emph{initial state} and \emph{accepting states}. Every edge is labeled by a letter of $\Sigma$. 
Recording the edge labels of any path from the initial state to some accepting state gives a word in $\Sigma^*$ that is \emph{accepted by the automaton}. A language $\cL \subseteq \Sigma^*$ is recognized by a finite automaton if it consists precisely of the words accepted by the automaton. 

\begin{ex}\label{ex:language poly ring} 
Let $\cL$ be the language on an alphabet $\Sigma=\{\tau,\alpha_1,\alpha_2\}$ with three letters defined as 
\begin{align*}
\cL & = \{ \tau^{k_1}\alpha_{i_1}\tau^{k_2}\dots \tau^{k_d}\alpha_{i_d}\tau^{k_{d+1}} \; \mid \; d \in \N_0,  i_1,\dots,i_d \in \{1,2\}, k_1,\dots, k_{d+1} \in \N_0, \\
& \hspace*{7.5cm} i_{j}\leq i_{j+1}\text{ if } k_j=0 \text{ with } 1 < j\leq d \}. 
\end{align*}
This is precisely the  language considered in \cite[Definition~4.2]{MN} with $c=(2)$ and $q=1$.  It is regular as it is recognized by the finite automaton in  \Cref{automaton poly ring}. 
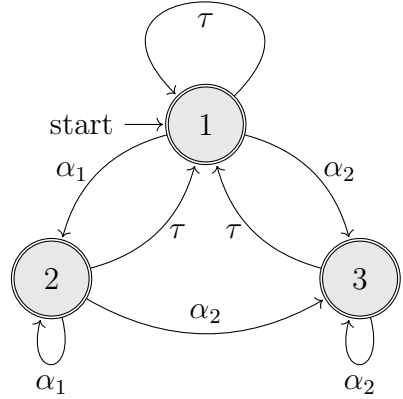
\begin{figure}[ht]
    \centering
\begin{tikzpicture}[shorten >=1pt,node distance=2.9cm,auto]
  \tikzstyle{every state}=[fill={rgb:black,1;white,10}]
   \node[state,initial,accepting]           (1)     {$1$};
   \node[state,accepting]           (3) [below right of=1]     {$3$};
  \node[state,accepting] (2) [below left of=1]  {$2$};
  \path[->]
    (1)   edge  [loop]         node {$\tau$} (1)
    (2)   edge  [loop below]         node {$\alpha_1$} (2)
    (3)   edge  [loop below ]         node {$\alpha_2$} (3)
   (1)   edge  [bend right]         node [left] {$\alpha_1$} (2)
    (2)   edge  [bend right]          node  {$\alpha_2$} (3)
     (2)   edge  [bend right]          node [right]{$\tau$} (1)
    (3)   edge  [bend left]          node[left] {$\tau$} (1)
     (1)   edge  [bend left]          node [right] {$\alpha_2$} (3);      
\end{tikzpicture}
\caption{The finite automaton for the language in \Cref{ex:language poly ring}.}  
\label{automaton poly ring}
\end{figure}
\end{ex}

A  \emph{weight function} on a language $\cL$ is a 
monoid homomorphism $ \rho \colon \cL \rightarrow \Mon (T)$, where $\Mon (T)$ denotes the set of monomials of a polynomial ring $T$ in finitely many variables $s_1,\ldots,s_k$. 
Its generating function is a formal power series 
\[
P_{\cL, \rho} (s_1,\ldots,s_k) = \sum_{w \in \cL} \rho (w).  
\]
If $\cL$ is a regular language it is a rational function (see, e.g., \cite{H} or \cite[Theorem 4.7.2]{St}). In fact, if $\cL$ is recognized by a finite automaton, then it can be computed explicitly as 
\begin{equation}
     \label{eq:automaton-HS}
P_{\cL, \rho} (s_1,\ldots,s_k) = \mathbf{u}^{T}(I_r-\sum_{a\in \Sigma}\rho(a)M_a)^{-1}\mathbf{e}_1, 
\end{equation}
where $I_r$ is the identity $r \times r$ matrix and $r$ is the number of states of the automaton, $\mathbf{u}$ is the indicator vector for the accepting states,  $\mathbf{e}_1$ is the indicator vector for the starting state and $M_a$ is the $0- 1$matrix to letter $a$ whose $(i,j)$ entry  is $1$ if there is an edge labeled $a$ from state $j$ to state $i$.

\begin{ex} 
      \label{ex:equivH poly ring} %[The equivariant Hilbert series for $\K[X_{[c]\times \N}]$.]
Consider the language $\cL$ of \Cref{ex:language poly ring} with  the weight function $\rho \colon \Sigma^*\rightarrow \Mon(\K[t,s])$, defined by $\rho(\alpha_1)=\rho(\alpha_2)=t$ and $\rho(\tau)=s$. Using the automaton in Figure~\ref{automaton poly ring}, we obtain for its generating function: 
\[
P_{\cL,\rho}(t,s)= \begin{bmatrix}
1&1&1
\end{bmatrix} \begin{bmatrix}
1-s & -s& -s \\ -t &1-t &0\\ -t& -t & 1-t
\end{bmatrix}^{-1}  \begin{bmatrix}
1\\0\\0
\end{bmatrix} =\dfrac{1}{(1-t)^2-s}.
\]
\end{ex}

Later, we will use regular languages to prove rationality of equivariant Hilbert series.

\section{Segre Products of Languages and Algebras} 
   \label{sec:Segre Products of Languages and Algebras}

We introduce the Segre product of formal languages. It is modeled after the Segre product of algebras. We consider the question whether the Segre product of two regular languages is again a regular language. 

Throughout this section we consider any two languages $\cL_{\sA} \subseteq \Sigma_A^*$ and 
$\cL_{\sB}  \subseteq \Sigma_B^*$ on alphabets $\Sigma_A = \{\tau_{1,1},\dots, \tau_{1,a} \alpha_1,\ldots,\alpha_p\}$ and 
$\Sigma_B = \{\tau_{2,1},\dots, \tau_{2,b},\beta_1,\ldots, \beta_q\}$ with $p+a$ and $q+b$ letters, respectively, where each 
alphabet is partitioned into two groups of letters.  We use different letters to distinguish the groups. 
Note that every word in 
$\Sigma_A^*$ can be written as  
\begin{align}
  \label{eq:words}
{\btau_1}^{\bk_1} \alpha_{i_1} {\btau_1}^{\bk_2} \alpha_{i_2} \ldots {\btau_1}^{\bk_d} \alpha_{i_d} {\btau_1}^{\bk_{d+1}}
\end{align}
with integers  $d \ge 0$,  $1 \le i_1,\ldots, i_d \le p$, and ${\btau_1}^{\bk_l}$ is some string that uses only letters $\tau_{1,1},\dots,\tau_{1,a}$. %\epsilon$, where $\epsilon$ is the empty word. 
Words in $\cL_A$ are words of the form  \eqref{eq:words} with conditions on $\bk_1,\dots,\bk_{d+1}$ and $i_1,\dots, i_d$; see for instance Examples  \ref{exa:easy language} and \ref{ex:language poly ring}. Similarly for $\Sigma_B^*$ and $\cL_B$.

\begin{defn} 
   \label{def:Segre language}
The \emph{Segre product} of $\cL_{\sA}$ and $\cL_{\sB}$ is the language  $\cL_{\sA}\boxtimes\cL_{\sB}$ on the alphabet 

\[
 \Sigma_A \boxtimes \Sigma_B =\{\tau_{1,1},\ldots,\tau_{1,a},\tau_{2,1},\dots,\tau_{2,b}, \gamma_{i,j} \; \mid \;  i \in [p], j \in [q] \}
 \] 
 with $p q +a+b$ letters defined by 
 \begin{align*}
 \hspace{2em}&\hspace{-1em}
 \cL_{\sA} \boxtimes \cL_{\sB}= \\
 & \{{\btau_1}^{\bk_1} {\btau_2}^{\bl_1} \gamma_{i_1,j_1} \ldots {\btau_1}^{\bk_d} {\btau_2}^{\bl_d}\:\gamma_{i_d, j_d}\:{\btau_1}^{\bk_{d+1}} {\btau_2}^{\bl_{d+1}} \; \mid \;  
{\btau_1}^{\bk_1} \alpha_{i_1} {\btau_1}^{\bk_2} \alpha_{i_2} \ldots \btau_1^{\bk_d} \alpha_{i_d} {\btau_1}^{\bk_{d+1}} \in \cL_{\sA} \\
& \hspace*{5cm} \text{ and } {\btau_2}^{\bl_1} \beta_{j_1} {\btau_2}^{\bl_2} \beta_{j_2} \ldots {\btau_2}^{\bl_d} \beta_{j_d} {\btau_2}^{\bl_{d+1}} \in \cL_{\sB} \}. 
 \end{align*}
\end{defn}

Observe that the Segre product $\cL_{\sA}\boxtimes\cL_{\sB}$  depends on the choice of partitions of the alphabets $ \Sigma_A$, $ \Sigma_B$. We suppress this dependency for ease of notation. 

Analogously, one can define a Segre product of more than two languages. We leave the details to the interested reader. 

\begin{ex}
The Segre product of $\{\tau_1,\alpha\}^*$ and  $\{\tau_2,\beta\}^*$ is
\begin{align*}
 \{\tau_1,\alpha\}^* \boxtimes \{\tau_2,\beta\}^* & = \{ \tau_1^{k_1} \tau_2^{l_1} \gamma \ldots \tau_1^{k_d}\tau_2^{l_d}\:\gamma\:\tau_1^{k_{d+1}} \tau_2^{l_{d+1}} \; \mid \;  k_{\nu}, l_{\nu} \in \N_0, \nu\in [d+1]\} \\
& = (\{\tau_1\}^* \{ \tau_2\}^* \gamma)^*. 
\end{align*} 
 Hence, this language is regular. 
\end{ex}

Some further instances in which the Segre product of two regular languages is regular have been established in \cite{MN}. In fact, this is always true as we show now. We thank Dietrich Kuske for the idea to prove the following result.

\begin{thm}
   \label{thm:regularity of Segre language}
If  $\cL_{\sA}$ and $\cL_{\sB}$ are regular languages then their Segre product $\cL_{\sA} \boxtimes \cL_{\sB}$ is a regular language as well. 
\end{thm}

\begin{proof}
Using the above notation, define a  monoid homomorphism $f_A \colon ( \Sigma_A \boxtimes \Sigma_B)^* \to \Sigma_A^*$ by $\tau_{1, j} \mapsto \tau_{1, j}, \ \tau_{2, j} \mapsto \epsilon$ and $\gamma_{i, j} \mapsto \alpha_i$, where $\epsilon$ denotes the empty word. Similarly, let  $f_B \colon ( \Sigma_A \boxtimes \Sigma_B)^* \to \Sigma_B^*$ be the  monoid homomorphism with 
$\tau_{1, j} \mapsto \epsilon, \ \tau_{2, j} \mapsto \tau_{2, j}$ and $\gamma_{i, j} \mapsto \beta_i$. Furthermore, consider the following language on $\Sigma_A \boxtimes \Sigma_B$, 
\[
\cL = \big (\{\tau_{1,1},\ldots, \tau_{1,a} \}^* \{\tau_{2,1},\ldots, \tau_{2,b} \}^* \{\gamma_{i, j} \; \mid \;  i \in [p], j \in [q]\} \big )^* \{\tau_{1,1},\dots, \tau_{1,a} \}^* \{\tau_{2,1},\ldots, \tau_{2,b} \}^*. 
\]
By the definition of the Segre product, one has 
\[
 \cL_{\sA} \boxtimes \cL_{\sB} = f_A^{-1} (\cL_{\sA}) \cap f_{\sB}^{-1} (\cL_B) \cap \cL. 
\]
Note that $\cL$ is a regular language. Since $\cL_{\sA}$ and $\cL_{\sB}$ are regular languages by assumption, \cite[Theorem 4.16]{HU} gives that $f_A^{-1} (\cL_{\sA})$ and $f_B^{-1} (\cL_{\sB})$ are also regular. Hence, as the intersection of regular languages, the language $\cL_{\sA} \boxtimes \cL_{\sB}$ is regular. 
\end{proof}

Using certain statistics on words, we now discuss some further properties of the Segre product. 
For any integer $d \ge 0$, define $\cL_{\sA}^d$, $\cL_{\sB}^d$ and $(\cL_{\sA} \boxtimes \cL_{\sB})^d$ as the set of words in $\cL_{\sA}$, $\cL_{\sB}$ and $\cL_{\sA} \boxtimes \cL_{\sB}$ with exactly $d$ $\alpha$-letters, $d$ $\beta$-letters and $d$ $\gamma$-letters, respectively. Moreover, let $(\cL_{\sA})_m$ be the set of words in $\cL_{\sA}$ with exactly $m$ letters $\tau_{1, i}$, $i \in [a]$. Define $(\cL_{\sB})_n$ by counting the occurrences of letters $\tau_{2, j}$ with $j \in [b]$. Similarly, for $(m, n) \in \N_0^2$, let 
$(\cL_{\sA} \boxtimes \cL_{\sB})_{m, n}$ be the set of words in $\cL_{\sA} \boxtimes \cL_{\sB}$ in which  letters $\tau_{1, i}$ occur exactly $m$ times and letters $\tau_{2, j}$ occur exactly $n$ times. Note that 
\[
(\cL_{\sA} \boxtimes \cL_{\sB})_{m, n} = \bigcup_{d \ge 0} (\cL_{\sA} \boxtimes \cL_{\sB})_{m, n}^d. 
\]
Analogous equalities are true for $\cL_{\sA}^d$, $\cL_{\sB}^d$. The definition of the Segre product of languages immediately implies the following observation. 

\begin{lm}
    \label{lem:bijection on Segre language} 
For any $m, n, d \in \N_0$, the map 
\[ 
(\cL_{\sA})^d_m\times(\cL_{\sB})^d_n \to (\cL_{\sA}\boxtimes \cL_{\sB})_{m,n}^d, \; (w_A, w_B) \mapsto w_{A, B}  
\]  
with 
\begin{align*}
w_A &=   \btau_1^{\bk_1}\:\alpha_{i_1}\:\dots \btau_1^{\bk_d}\:\alpha_{i_d}\:\btau_1^{\bk_{d+1}},  \\
w_B &=  \btau_2^{\bl_1}\:\beta_{j_1}\:\dots \btau_2^{\bl_d}\:\beta_{j_d}\:\btau_2^{\bl_{d+1}}, \; \text{ and } \\
w_{A, B} & =  \btau_1^{\bk_1}\btau_2^{\bl_1}\:\gamma_{i_1, j_1}\:\dots \btau_1^{\bk_d}\btau_2^{\bl_d}\:\gamma_{i_d, j_d}\:\btau_1^{\bk_{d+1}}\btau_2^{\bl_{d+1}}
\end{align*}
 is bijective. 
 \end{lm}

Consider weight functions $\rho_A \colon \Sigma_A^* \to \Mon (\K[s_1, t])$ and 
$\rho_B \colon \Sigma_B^* \to \Mon (\K[s_2, t])$, defined by $\rho_A (\tau_{1, i}) = s_1$ for $i \in [a]$, $\rho_A (\alpha_k ) = t$ for $k \in [p]$ and $\rho_B (\tau_{2,j}) = s_2$ for $j \in [b]$, $\rho_B (\beta_k ) = t$ for $k \in [q]$. These weight functions induce weight functions on $\cL_{\sA}$ and  $\cL_{\sB}$, respectively. Recall that the resulting generating function for $\cL_{\sA}$ is 
the \emph{formal power series} 
\[ 
P_{\cL_{\sA},\rho_{A}} (s_1, t) =  \sum_{w\in \cL_{\sA}} \rho_A (w) \in \mathbb{Z}[[s_1,t]]
\]
and similarly $P_{\cL_{\sB},\rho_{B}} (s_2, t) \in \mathbb{Z}[[s_2,t]]$. 

Define a weight function $\rho_{A, B} \colon  (\Sigma_A \boxtimes \Sigma_B)^* \to \Mon (\K [s_1, s_2, t])$ by 
\[
\rho_{A, B} (\tau_{1, i}) = s_1 \text{ if } i \in [a], \ \rho_{A, B} (\tau_{2, j}) = s_2 \text{ if } j \in [b], \text{ and } \rho_{A, B} (\gamma_{k, l}) = t \text{ if } k \in [p], l \in [q]. 
\]
Its generating function is determined by those of $\rho_A$ and $\rho_B$. More precisely, one has the following equality. 

\begin{prop} 
    \label{prop:Hilb Segre formula} 
\begin{align*}
 {P}_{\cL_{\sA}\boxtimes\cL_{\sB},\rho_{A, B}} (s_1,s_2, t) &=  \sum_{m,n,d \in \N_0^3 } \left (\sum_{w \in (\cL_{\sA})_{m}^d} \rho_A (w) \right ) \cdot \left (\sum_{w \in (\cL_{\sB})_{n}^d} \rho_B (w) \right ) \cdot  t^{-d}
\end{align*}
\end{prop}

\begin{proof}
For any $w\in(\cL_{\sA}\boxtimes\cL_{\sB})_{m,n}^d$, we have $\rho_{A, B} (w)=t^ds_1^ms_2^n$. Using also \Cref{lem:bijection on Segre language}, it follows 
\begin{align*}
 {P}_{\cL_{\sA}\boxtimes\cL_{\sB},\rho_{A, B}} (s_1,s_2, t)  &= \sum_{w\in \cL_{\sA}\boxtimes\cL_{\sB}} \rho(w)  
 =  \sum_{m,n,d \in \N_0^3} \sum_{w\in (\cL_{\sA}\boxtimes\cL_{\sB})_{m,n}^d} \rho_{A, B} (w) \\
 &=\sum_{m,n,d \in \N_0^3 } \#(\cL_{\sA}\boxtimes\cL_{\sB})_{m,n}^d \cdot  t^ds_1^ms_2^n\\
 &=\sum_{m,n,d \in \N_0^3 } \#(\cL_{\sA})_{m}^d \cdot \#(\cL_{\sB})_{n}^d \cdot  t^ds_1^ms_2^n\\
&=\sum_{m,n,d \in \N_0^3 } \left (\sum_{w \in (\cL_{\sA})_{m}^d} \rho_A (w) \right ) \cdot \left (\sum_{w \in (\cL_{\sB})_{n}^d} \rho_B (w) \right ) \cdot  t^{-d}, 
\end{align*}
as claimed. 
\end{proof}

We now relate the above construction to equivariant Hilbert series. We begin by recalling the Segre product of algebras. 

Temporarily using new notation, let $ A = \K[a_1,\ldots,a_s] \subset R$  and $B = \K[b_1,\ldots,b_t] \subset S$ be subalgebras of standard graded polynomial rings  $R = \K[x_1,\ldots,x_m]$ and  $S = \K[y_1,\ldots,y_n]$ that are generated by monomials $a_1,\ldots,a_s$ of degree $d_1$ and  monomials $b_1,\ldots,b_t$ of degree $d_2$, respectively. 
Denote by $C$ the subalgebra of $\K[x_1,\ldots,x_m, y_1,\ldots,y_n]$ that is generated by all monomials $a_i b_j$ with $i \in [s]$ and $j \in [t]$. Using the gradings induced from the corresponding polynomials rings one has (see, e.g., \cite[Lemma 5.1]{MN}), for all $k \in \Z$, 
\[
\dim_{\K} [C]_{k (d_1 + d_2)} = \dim_{\K} [A]_{k d_1} \cdot \dim_{\K} [B]_{k d_2}. 
\]
Regrading $A$ and $B$ as standard graded algebras that are generated in degree one by setting $[A]_j = A \cap [R]_{j d_1}$ and $[B]_j = B \cap [S]_{j d_2}$, their Segre product becomes an algebra generated in degree one with grading 
$[C]_j = C \cap [R \otimes_K S]_{j (d_1 +d_2)}$. We denote $C$ with this grading by $A \boxtimes B$. Observe that in the new grading, the last equality becomes 
\begin{equation}
    \label{eq:hilb Segre}
\dim_{\K} [A \boxtimes B]_{k} = \dim_{\K} [A]_{k} \cdot \dim_{\K} [B]_{k}, 
\end{equation}
which justifies to call $A \boxtimes B$ the \emph{Segre product} of  $A$ and $B$. 

We now return to the notation used in \Cref{sec:prep} though we use a more general setup. Fix constants $m_0 \in \Z$ and $c_1 \in \N_0$ and denote by  $\sA=(A_m)_{m \ge m_0}$  a family of  algebras, where every $A_m$
is a finitely generated monomial subalgebra of the polynomial ring $R_m=\K[x_{i, j} \; \mid \; i \in [c_1], j \in [m]]$ and its generating monomials all have the same degree, say $d_1$, considered as monomials in $R = {\displaystyle \lim_{\longrightarrow}}\, R_m = \K[x_{i, j} \; \mid \; i \in [c_1], j \in \N]$. If, for any $n \ge m \ge m_0$,  there are $\K$-algebra homomorphisms $f_{m, n} \colon A_m \to A_n$ that form a direct system we 
we call $\sA=(A_m)_{m \ge m_0}$ a \emph{1-filtration} or simply a \emph{filtration} of algebras. 

We formalize the relation of such a family to a language on the alphabet $\Sigma_A = \{\tau_{1,1},\ldots,\tau_{1, a}, \alpha_1,\ldots,\alpha_p\}$. 

\begin{defn}
  \label{def:represented by language} 
We say that a 1-filtration $\sA$ is \emph{represented by a language} $\cL_{\sA} \subset \Sigma_A^*$, if there is an integer $\tilde{m}$ such that, for any $m \ge \tilde{m}$ and any $d \in \N_0$, there is a bijection 
\[
(\bm_A)_m^d \colon [\Mon (A_m)]_d \to (\cL_{\sA})_m^d, 
\]
where $[\Mon (A_m)]_d$ denotes the set of monomials in $A_m$ of (internal) degree $d$, that is, of degree $d d_1$ when considered as elements of $R$. 
\end{defn}

For fixed integers $n_0$ and $c_2 \ge 1$, consider another  family of  algebras $\sB=(B_n)_{n \ge n_0}$  a, where every $B_n$
is a finitely generated monomial subalgebra of the polynomial ring $S_n=\K[y_{i, j} \; \mid \; i \in [c_2], j \in [n]]$ and its generating monomials all have the same degree, say $d_2$, considered as monomials in $S = {\displaystyle \lim_{\longrightarrow}}\, S_n = \K[x_{i, j} \; \mid \; i \in [c_2], j \in \N]$. 

It is worth extending the concept of a 1-filtration to an $r$-filtration. We will work this out explicitly only  in the case $r=2$ and leave the more general case to the reader. A family of algebras $\sC = (C_{m, n})_{m, n}$ is called a \emph{2-filtration} if $(C_{m, n})_{m}$ and $(C_{m, n})_{n}$ are 1-filtrations. 

Consider a language $\cL$ on an alphabet $\Sigma$ containing at least the letters $\tau_{1,1},\ldots, \tau_{1,a}$ and  $\tau_{2,1},\ldots, \tau_{2,b}$. Denote by $\cL_{m, n}^d$ the set of words in $\cL$ with exactly $m$ occurrences of letters $\tau_{1, i}$ with $i \in [a]$, $n$ occurrences of letters $\tau_{2, j}$ with $j \in [b]$ and $d$ letters other than $\tau_1$ or $\tau_2$. We say that a 2-filtration $\sC$ is \emph{represented by} $\cL$, if,  for any $m \gg 0$, $n \gg 0$  and any $d \in \N_0$, there is a bijection 
\[
 [\Mon (C_{m, n})]_d \to \cL_{m, n}^d.  
\]

The following result justifies calling the language $\cL_{\sA}\boxtimes\cL_{\sB}$ a Segre product. 

\begin{thm}
    \label{prop:relate Segre products}
If 1-filtrations of algebras $\sA=(A_m)_{m \ge m_0}$ and $\sB=(B_n)_{n \ge n_0}$ are represented by languages  $\cL_{\sA}$ and $\cL_{\sB}$, respectively, then one has: 
\begin{itemize}

\item[(a)]  
For every $m, n \gg 0$ and any $d \in \N_0$, there are bijections
\[
[\Mon (A_m \boxtimes B_n)]_d \to (\cL_{\sA}\boxtimes\cL_{\sB})_{m, n}^d,  
\]
that is, 
the 2-filtration $(A_m \boxtimes B_n)_{m \ge m_0, n \ge n_0}$ is represented by $\cL_{\sA}\boxtimes\cL_{\sB}$. 

\item[(b)] If $\cL_{\sA}$  and $\cL_{\sB}$ are regular languages, then the equivariant Hilbert series of $\sA \boxtimes \sB$, that is, 
\begin{align*}
{\rm equivH}_{\sA \boxtimes \sB} (s_1,s_2, t) & = \sum_{m \ge m_0, n \ge n_0} H_{A_m\boxtimes B_n} (t) \cdot s_1^m s_2^n \\
& = 
\sum_{m \ge m_0, n \ge n_0} \sum_{d \ge 0} \dim_{\K} [A_m \boxtimes B_n]_d \cdot t^d s_1^m s_2^n 
\end{align*}
is rational. 
\end{itemize}
\end{thm}

\begin{proof} 
(a) 
Since every algebra $A_m$, $B_n$ is monomial, any degree component of any of the involved algebras has a $\K$-vector space basis consisting of monomials. Thus, Equality \eqref{eq:hilb Segre} gives bijections 
\[
[\Mon (A_m \boxtimes B_n)]_d \to  [\Mon (A_m)]_d \times [\Mon (B_n)]_d. 
\]
By assumption, for any $m \gg 0$ and $d \in \N_0$, there is a bijection $(\bm_A)_m^d \colon [\Mon (A_m)]_d \to (\cL_{\sA})_m^d$ and, similarly, there is a bijection $(\bm_B)_n^d \colon [\Mon (B_n)]_d \to (\cL_{\sB})_n^d$ whenever $n \gg 0$ and $d \in \N_0$. Moreover, by \Cref{lem:bijection on Segre language}, there are bijections 
$(\cL_{\sA})^d_m\times(\cL_{\sB})^d_n \to (\cL_{\sA}\boxtimes \cL_{\sB})_{m,n}^d$. Now the claim follows. 

(b)  Part (a) gives
\begin{align*}
{\rm equivH}_{\sA \boxtimes \sB} (s_1,s_2, t) 
& = \sum_{m \ge m_0, n \ge n_0} \sum_{d \ge 0} \dim_{\K} [A_m \boxtimes B_n]_d \cdot t^d s_1^m s_2^n \\
& = \sum_{m \ge m_0, n \ge n_0} \sum_{d \ge 0} \# [\Mon (A_m \boxtimes B_n)]_d \cdot t^d s_1^m s_2^n \\
& = \sum_{m \ge m_0, n \ge n_0} \sum_{d \ge 0} \# (\cL_{\sA} \boxtimes \cL_{\sB})_{m, n}^d \cdot t^d s_1^m s_2^n \\ 
\end{align*}
Applying \Cref{lem:bijection on Segre language},  we get. 
\begin{align*}
{\rm equivH}_{\sA \boxtimes \sB} (s_1,s_2, t) 
& =  \sum_{m \ge m_0, n \ge n_0} \sum_{d \ge 0} \# (\cL_{\sA})_{m}^d \cdot \# (\cL_{\sB})_{n}^d \cdot t^d s_1^m s_2^n.  \end{align*}
Hence \Cref{prop:Hilb Segre formula} implies
\begin{align*}
{\rm equivH}_{\sA \boxtimes \sB} (s_1,s_2, t) 
& =  \sum_{m \ge m_0, n \ge n_0} \sum_{d \ge 0} \# (\cL_{\sA})_{m}^d \cdot \# (\cL_{\sB})_{n}^d \cdot t^d s_1^m s_2^n  \\
%&=\sum_{m,n,d \in \N_0^3 } \#(\cL_{\sA})_{m}^d \cdot \#(\cL_{\sB})_{n}^d \cdot  t^ds_1^ms_2^n\\
&=\sum_{m \ge m_0, n \ge n_0} \sum_{d \ge 0} \left (\sum_{w \in (\cL_{\sA})_{m}^d} \rho_A (w) \right ) \cdot \left (\sum_{w \in (\cL_{\sB})_{n}^d} \rho_B (w) \right ) \cdot  t^{-d} \\
& = {P}_{\cL_{\sA}\boxtimes\cL_{\sB},\rho_{A, B}} (s_1,s_2, t). 
\end{align*}
Since $\cL_{\sA}\boxtimes\cL_{\sB}$ is a regular language by \Cref{thm:regularity of Segre language}, ${P}_{\cL_{\sA}\boxtimes\cL_{\sB},\rho_{A, B}} (s_1,s_2, t)$ is a rational function. 
\end{proof} 

We close this section with an application to the filtration of Segre products of polynomial rings. 

\begin{ex} 
      \label{exa:equivH poly ring}
(i)  
First we consider the 1-filtration $\sA=(A_n)_{n\in \N}$ with $A_n=\K[X_{[c] \times [n]}]$, where $c \ge 1$ is any fixed integer.  
Generalizing \Cref{ex:language poly ring}, $\sA$ is represented by the language  $\cL$ on an alphabet $\Sigma=\{\tau,\alpha_1,\ldots,\alpha_c\}$  defined as 
\begin{align*}
\cL & = \{ \tau^{k_1}\alpha_{i_1}\tau^{k_2}\dots \tau^{k_d}\alpha_{i_d}\tau^{k_{d+1}} \; \mid \; d \in \N_0,  i_1,\dots,i_d \in [c], k_1,\dots, k_{d+1} \in \N_0, \\
& \hspace*{7.5cm} i_{j}\leq i_{j+1}\text{ if } k_j=0 \text{ with } 1 < j\leq d \}. 
\end{align*}
Indeed, if $\cL_{n,d}$ denotes the set of words in $\cL$ that use exactly $n$ letters $\tau$ and $d$ of letters $\alpha_1,\dots,\alpha_c$, then one shows that the map $\bm \colon \cL \rightarrow \bigcup_{n \ge 1, d \ge 0} [\Mon(A_n)]_d$ with 
$\bm (\tau)=1$ and 
$\bm (\tau^k\alpha_i)=x_{i,k+1}$ for $k\in \N_0$ and $i\in [c]$ induces bijections $\cL_{n-1}^d\rightarrow [\Mon(A_n)]_d$.  

(ii) 
Second, fix integers $c_1, c_2 \ge 1$ and consider the 2-filtration $\sC = (C_{m, n})_{m, n \in \N}$ with 
\[
C_{m, n} = \K[X_{[c_1] \times [m]}] \boxtimes \K[Y_{[c_2] \times [n]}]. 
\]
Let $\cL_1$ and  $\cL_2$ be the languages representing the 1-filtrations  $(\K[X_{[c_1]\times [m]}])_{m \in \N}$ and 
$(\K[Y_{[c_2]\times [n]}])_{n \in \N}$, respectively,  as in (i). Their Segre product $\cL = \cL_1 \boxtimes \cL_2$ is given explicitly as 
\begin{align*}
 \cL=\{&\tau_1^{k_{1,1}} \tau_2^{k_{2,1}}\gamma_{i_{1,1}, i_{2,1}} \dots \dots  \tau_1^{k_{1,d}} \tau_2^{k_{2,d}}\gamma_{i_{1,d},i_{2,d}}\tau_1^{k_{1,d+1}} \tau_2^{k_{2,d+1}} \mid \\ 
 & \hspace*{3.5cm} \tau_j^{k_{j,1}}\alpha_{i_{1,j}}\dots \dots  \tau_j^{k_{j,d}}\alpha_{i_{j,d}}\tau_j^{k_{j,d+1}}\in \cL_j \text{ for } j\in \{1, 2\}\}
\end{align*}
It is a regular language representing $\sC$. Indeed, the map $\bm \colon \cL \rightarrow \bigcup_{m, n \ge 1, d \ge 0} [\Mon(C_{m,n})]_d$ with 
$\bm (\tau)=1$ and 
$\bm (\tau_1^{k_1}  \tau_2^{k_2} \gamma_{j_1, j_2})=x_{j_1,k_1+1} y_{j_2, k_2 +1}$ induces bijections $\cL_{m-1, n-1}^d\rightarrow [\Mon(C_{m,n})]_d$. Using the weight function $\rho \colon \Sigma^*\rightarrow \K[s_1,s_2, t]$ with  
$\rho(\gamma_{i_1, i_2})=t$ and $\rho (\tau_j) = s_j$, one obtains that $\equivH_{\sA}(s_1,s_2,t)=s_1 s_2 \cdot P_{\cL,\rho}(s_1,s_2,t)$
is a rational function in $\mathbb{Q}(s_1,s_2,t)$. For details we refer the reader to \cite{MN}.
\end{ex}

%%%%%%%%%%%%%%%%%

\section{Tensor products} 
   \label{sec:tensor product} 

Locally, a Segre product is given as a tensor product of algebras. We complement the previous section by 
considering 2-families defined by tensor products of 1-families.  We show that they have a rational equivariant Hilbert series if the  factors are represented by regular languages. 

We continue to use the notation introduced in \Cref{sec:Segre Products of Languages and Algebras}. In particular, $\cL_{\sA} \subseteq \Sigma_A^*$ and 
$\cL_{\sB}  \subseteq \Sigma_B^*$ are languages on alphabets $\Sigma_A = \{\tau_{1,1},\dots, \tau_{1,a} \alpha_1,\ldots,\alpha_p\}$ and 
$\Sigma_B = \{\tau_{2,1},\dots, \tau_{2,b},\beta_1,\ldots, \beta_q\}$.  We partition the concatenation $\cL_{\sA} \cL_{\sB}$ into subsets 
$(\cL_{\sA} \cL_{\sB})_{m, n}^d$, where $(\cL_{\sA} \cL_{\sB})_{m, n}^d$ consists of words with exactly $m$ occurrences of the letters $\tau_{1,1},\ldots,\tau_{1,a}$, $n$ occurrences of $\tau_{2,1},\ldots, \tau_{2, b}$ and $d$ letters other than $\tau_{1,i}$ or $\tau_{2,j}$.

The definitions imply the following observation. We use $M \uplus N$ to denote the disjoint union of sets $M$ and $N$. 

\begin{lm}
    \label{lem:bijection on concatenation} 
For any $m, n, d \in \N_0$, the map 
\[ 
\biguplus_{d_1 + d_2 = d} (\cL_{\sA})^{d_1}_m \times (\cL_{\sB})^{d_2}_n \to (\cL_{\sA} \cL_{\sB})_{m,n}^d, \; (w_A, w_B) \mapsto w_{A} w_B,   
\]  
is bijective.
\end{lm}

The main result of this section is analogous to \Cref{prop:relate Segre products} and \Cref{thm:Segre product with rat hilb}.    

\begin{prop}
    \label{prop:relate tensor products}
Given 1-filtrations of monomial algebras $\sA=(A_m)_{m \ge m_0}$ and $\sB=(B_n)_{n \ge n_0}$, one has: 
\begin{itemize}

\item[(a)] If $\sA$ and $\sB$ are represented by languages  $\cL_{\sA}$ and $\cL_{\sB}$, respectively, then the 2-filtration 
$\sA \otimes \sB = (A_m \otimes_{\K} B_n)_{m \ge m_0, n \ge n_0}$ is represented by $\cL_{\sA} \cL_{\sB}$. 

\item[(b)] If $\cL_{\sA}$ and $\cL_{\sB}$ are regular languages then the equivariant Hilbert series 
\begin{align*}
{\rm equivH}_{\sA \otimes \sB} (s_1,s_2, t) & = \sum_{m \ge m_0, n \ge n_0} H_{A_m\otimes B_n} (t) \cdot s_1^m s_2^n \\
& = 
\sum_{m \ge m_0, n \ge n_0} \sum_{d \ge 0} \dim_{\K} [A_m \otimes B_n]_d \cdot t^d s_1^m s_2^n 
\end{align*}
is rational.
\end{itemize}
\end{prop} 

\begin{proof}
The argument for Claim (a) is analogous to the proof of \Cref{prop:relate Segre products}. The only difference is that this time we use \Cref{lem:bijection on concatenation} and bijections
\[
[\Mon (A_m \otimes B_n)]_d \to  \bigcup_{d_1 + d_2 = d}  [\Mon (A_m)]_{d_1} \times [\Mon (B_n)]_{d_2}. 
\]

For Claim (b), we define a weight function $\rho \colon \Sigma^*_{\sA}\otimes\Sigma^*_{\sB}\rightarrow \Z[t,s_1,s_2]$ by 
\[
\rho(\alpha_{1})=\dots=\rho(\alpha_{p})=\rho(\beta_{1})=\dots=\rho(\beta_{q})=t, \text{ and }
\rho(\tau_{1,i})=s_1,\; \rho(\tau_{2,j})=s_2 \text{ for } i\in [a],j\in [b].
\]
It follows 
\begin{align*}
{\rm equivH}_{\sA \otimes \sB} (s_1,s_2, t) 
& = \sum_{m \ge m_0, n \ge n_0} \sum_{d \ge 0} \dim_{\K} [A_m \otimes B_n]_d \cdot t^d s_1^m s_2^n \\
& = \sum_{m \ge m_0, n \ge n_0} \sum_{d \ge 0} \# [\Mon (A_m \otimes B_n)]_d \cdot t^d s_1^m s_2^n \\
& = \sum_{m \ge m_0, n \ge n_0} \sum_{d \ge 0} \# (\cL_{\sA} \cL_{\sB})_{m, n}^d \cdot t^d s_1^m s_2^n \\
& = \sum_{m \ge m_0, n \ge n_0} \sum_{d \ge 0} \sum_{w\in (\cL_{\sA} \cL_{\sB})_{m,n}^d} \rho (w)\\ 
& = {P}_{\cL_{\sA} \cL_{\sB},\rho} (s_1,s_2, t).  
\end{align*}
As concatenation of regular languages, $ \cL_{\sA} \cL_{\sB}$ is regular. Thus, ${P}_{\cL_{\sA} \cL_{\sB},\rho} (s_1,s_2, t)$ is a rational function (see \Cref{eq:automaton-HS}). 
%Thus, we conclude by \Cref{thm:reg has rat}. 
\end{proof}

\begin{ex}
For any $(p, q) \in \N_0$, the 2-filtration $(C_{m, n})_{m,n \in \N^2}$ with 
\[ 
C_{m, n} = \K[x_i^2,x_ix_{i+1},\dots,x_ix_{i+p}, y_j^2,y_j y_{j+1},\dots,y_j u_{j+q} \:\mid\: i \in  [m], j \in [n] ].
\]
has a rational equivariant Hilbert series. 
\end{ex}

%%%%%%%%%%%%%%%%%%%%%%%%%%%%%%

\section{Equivariant Hilbert Series of Some Algebras}
\label{sec:1}

We consider an infinite family of filtrations of monomial algebras. Each filtration has a colimit that is finitely generated up to shifting as an algebra. We show that each filtration in the family has a rational equivariant Hilbert series and that such rationality  is also true for the filtration obtained by taking Segre products of any two filtrations in the family. 

Given $c\in \N_0$,  denote by  is a subalgebra of the polynomial ring $R_n=\K[x_1,\dots,x_{n+c}]$. 
Using the natural inclusions $A_m \to A_n$ if $m \le n$, the colimit
\[
A = {\displaystyle \lim_{\longrightarrow}}\, A_n = \K[x_i x_j \; \mid \; 1 \le i \le j \le i+c ]
\]
is a subalgebra of the polynomial ring $R=\K[x_i  \; \mid \; i \in \N ]$ in infinitely many variables. Observe that $A$  is finitely generated by $G=\{x_1^2,x_1x_2,\ldots,x_1x_{1+c}\}$ up to  shifting.

Denote by $\Mon (R)$ the set of monomials in $R$. Each such monomial 
 can be uniquely written as a string of variables in the non-commutative polynomial ring $NR=\K \langle x_1,x_2, \dots \rangle $, where indices  appear in a non-decreasing order from left to wright. Denote the set of all such strings by $\mathcal{S}$. Thus, we get a bijective map 
 \[
 \mathbf{s} \colon \Mon(R) \to \mathcal{S}
 \]
 that maps a monomial  $m \in \Mon (R)$ onto its  \emph{string presentation}  $\mathbf{s}(m)$. 
 For example, the string presentation of   $x_1x_2\,x_1x_3\,x_0x_0\,x_3x_3$ is  $x_{\color{red}{0}}x_{\color{red}{0}}x_{\color{red}{1}}x_{\color{red}{1}}x_{\color{red}{2}}x_{\color{red}{3}}x_{\color{red}{3}}x_{\color{red}{3}} \in \mathcal{S}$. 

We know describe the string presentations of monomials in $A$. 

\begin{lm}
   \label{lem:string present}
One has 
\[
 \mathbf{s} (\Mon (A)) = \{ x_{i_1}x_{i_2}\dots x_{i_{2d-1}}x_{i_{2d}} \; \mid \; d \in \N_0,   i_1 \leq i_2\leq \dots \leq i_{2d}, 
 i_{2 k} - i_{2k -1} \le c \text{ if } 1 \le k \le d \}. 
\]
\end{lm}

\begin{proof} 
Denote the right-hand side in the claimed equality by $M$. 
Since the generators of $A$ up to a shift are $x_1^2,x_1x_2,\ldots,x_1x_{1+c}$, every monomial in $A$ can be written as 
\begin{align}
    \label{eq:first reordering} 
m = (x_{i_1}x_{i_1 +j_1}) \cdot (x_{i_2}x_{i_2 +j_2}) \cdots (x_{i_{2d-1}}x_{i_{2-1} + j_{2d-1}}) 
\end{align}
with $i_1 \leq i_2\leq \dots \leq i_{d}$ and $0 \le j_k \le c$ for $k \in [d]$. It follows that $M \subseteq \mathbf{s} (\Mon (A))$. 

In order to prove the reverse inclusion we use induction on $d \ge 0$ to show that any monomial $m$ as in \eqref{eq:first reordering} has a string presentation in $M$. If $d = 0$ or  $d = 1$ this is clear. 

Let $d \ge 2$. We consider two cases. 

\emph{Case 1: } Assume $i_1 + j_1 \le i_2$. Write $m = (x_{i_1}x_{i_1 +j_1}) \cdot m'$. Then 
$ \mathbf{s} (m) = x_{i_1}x_{i_1 +j_1} \mathbf{s} (m')$. Moreover, we have $\mathbf{s} (m') \in M$ by induction. Since $\mathbf{s} (m')$ begins with $x_{i_2}$, it follows that $x_{i_1}x_{i_1 +j_1} \mathbf{s} (m') =  \mathbf{s} (m)$ is in $M$, as desired. 

\emph{Case 2: } Assume $i_1 + j_1 > i_2$. Write $m = (x_{i_1}x_{i_2}) \cdot m'$, and so $m' = (x_{i_1 +j_1} x_{i_2 +j_2}) \cdot (x_{i_3}x_{i_3 +j_3}) \cdots (x_{i_{2d-1}}x_{i_{2-1} + j_{2d-1}})$. Since $i_1 \le i_2 < i_1 + j_1 \le i_1 + c$ we see that 
$\mathbf{s} (x_{i_1}x_{i_2}) = x_{i_1}x_{i_2}$ is in $M$. Furthermore, the assumptions imply that $| i_2 +j_2 - (i_1 +j_1)| \le c$. It follows that, possibly after permuting $x_{i_1 +j_1}$ and  $x_{i_2 +j_2}$, the monomial $m'$ is as in \eqref{eq:first reordering}. Hence the induction hypothesis gives $\mathbf{s} (m') \in M$, which yields $x_{i_1}x_{i_2} \mathbf{s} (m') = \mathbf{s} (m) \in M$, completing the argument.  
\end{proof}

Now we want to enumerate the monomials in $A$ using the words of a suitable formal language on the alphabet $\Sigma = \{\tau, \alpha_0,\ldots, \alpha_c\}$ with $c+2$ letters. 
Consider a shift operator $T:\textrm{Mon}(A) \rightarrow \textrm{Mon}(A)$ defined by 
\[T(x_{i})=x_{i+1},\] 
and extended multiplicatively to $\textrm{Mon}(A)$.  Define a map $\textbf{m} \colon \Sigma^{\star} \rightarrow \textrm{Mon}(A)$ recursively using the following three rules: \\
\indent (a) $\textbf{m} (\epsilon)=1$,  \quad
(b)  $\textbf{m} (\alpha_{i} w )=x_1 x_{i+1}\textbf{m}(w)$,  \quad
(c)  $\textbf{m} (\tau w)=T(\textbf{m}(w))$ \quad 
if $w\in \Sigma^*$, 
where $\epsilon$ denotes the empty word. 
For example, $\mathbf{m}(\alpha_i)=x_1x_{i+1}$,  $\mathbf{m}(\tau^k\alpha_i)=x_{k+1} x_{i+k+1}$ and \(\mathbf{m}(\alpha_1\tau\alpha_0\alpha_2\tau^2)=x_1 x_2^4  x_4\) if $c\geq 2$. 

If $c>0$, the map $\textbf{m}$ is not injective since the variables $x_i$ commute. For example, $\mathbf{m}(\alpha_i\alpha_j)=\mathbf{m}(\alpha_j\alpha_i)$. Thus, we consider the language $\cL \subset \Sigma^*$ defined as 
\begin{align}
    \label{eq:language for factor} 
\cL &= 
\{ \tau^{k_1}\alpha_{i_1}\tau^{k_2}\alpha_{i_2}\dots \tau^{k_d}\alpha_{i_d}\tau^{k_{d+1}} \; \mid \; k_1, k_{d+1} \in \N_0, 0 \le  i_1,\dots,i_d\leq c, k_{j+1}\geq i_{j} \text{ if }  1 \le j < d \} .
\end{align}

Let $\cL_n^d$ be the collection of words in $\cL$  that use exactly $d$ times one of the letters  $\alpha_0,\ldots,\alpha_c$, and $n$ times the letter  $\tau$. Denote by $[\Mon (A_n)]_d$  the set of  monomials in $A_n$ of degree $2d$ (any monomial in $A_n$ must be of even degree).

\begin{lm} 
   \label{lem:bijection} 
For every $d \in \N_0$ $n \in \N$, there is a bijection $[\Mon (A_n)]_d \to  \cL_{n-1}^d$. 
\end{lm}{}

\begin{proof}
Using the above bijection $\mathbf{s}$, it suffices to show that the map $\mathbf{w}_n^d \colon \mathcal{S}_{n}^d = \mathbf{s} ([\Mon (A_n)]_d) \to  \cL_{n-1}^d$ is bijective, where $\mathbf{w}_n^d$ maps a string 
\[
s=x_{i_1}x_{i_2}\dots x_{i_{2d-1}}x_{i_{2d}} \; \text{ with } 1 \le i_1 \leq i_2\leq \dots \leq i_{2d} \le n+c, 
 i_{2 k} - i_{2k -1} \le c \text{ if } 1 \le k \le d 
\] 
onto the word 
\[
\mathbf{w}_n^d (s)=\tau^{i_1 - 1}\alpha_{i_2-i_1}\tau^{i_3-i_1}\alpha_{i_4-i_3}\dots\tau^{i_{2d-1}-i_{2d-3}}\alpha_{i_{2d}-i_{2d-1}}\tau^{n-i_{2d-1}}.
\]
Notice that the right-hand side is indeed in  $\cL_n^d$. 

The map $\mathbf{m}\circ \mathbf{s} \colon \cL \to  \mathcal{S}$ assigns to any word 
 \[
 w=\tau^{k_1-1}\alpha_{i_1}\tau^{k_2}\alpha_{i_2}\dots \tau^{k_d}\alpha_{i_d}\tau^{k_{d+1}} \in \cL_{n-1}^d
 \]
the string 
\[
(\mathbf{m}\circ \mathbf{s}) (w)=x_{k_1}x_{i_1+k_1}x_{k_1+k_2}x_{i_2+k_1+k_2}\dots x_{k_1+..+k_d}x_{i_d+k_1+..+k_d}.
\]
Since $w \in \cL_{n-1}^d$  means $k_{j+1}\geq i_{j} \text{ if }  1 \le j < d$ and $k_1 + \cdots + k_{d+1} = n$,  we get 
\[
 k_1\leq i_1+k_1\leq k_1+k_2\leq \dots \leq k_1+ \cdots +k_d\leq i_d+k_1+ \cdots + k_d \le c+n, 
 \]
 which shows that $(\mathbf{m}\circ \mathbf{s}) (w)$ is in  $\mathcal{S}_{n}^d$.  

Finally, one checks that the restriction of  $(\mathbf{m}\circ \mathbf{s})$ to  $\mathcal{S}_{n}^d$ and the map $\mathbf{w}_n^d$ are inverse to each other, which establishes the desired bijection.  
\end{proof}

Recall that the \emph{equivariant Hilbert} series of  $\sA$ is 
\begin{align*}
{\rm{equivH}}_{\mathscr{A }} (s, t) & = \sum_{n \in \mathbb{N}} H_{A_{n}} (t) \cdot s^n, 
\end{align*}
where $H_{A_n}(t) = \sum_{d \ge 0} \dim_{\K} [A_n]_d t^d$ is the Hilbert series for $A_n$ and, usuing the fact that any polynomial in $A_n$ has even degree,  we grade  $A_n$ by 
$[A_n]_d = A_n \cap [R_n]_{2d}$. 

Combining the above preparations we can now show the  main result of this section.  

\begin{thm}
\label{thm:thm1}
The equivariant Hilbert series of $\mathscr{A}$ is rational.
\end{thm}{}

\begin{proof}
Consider the language $\cL$ given in \eqref{eq:language for factor} and define a weight function $\rho \colon \Sigma^*\rightarrow \Z[t,s]$ on the free monoid  $\Sigma^*$ by
\begin{align*}
   \rho(\alpha_0)=\dots=\rho(\alpha_c)=t \; \text{ and } \;  \rho(\tau)=s.
\end{align*}
Thus, for every $w\in\cL^d_n$, one has $\rho(w)=t^ds^n$. This implies for the generating function by using also \Cref{lem:bijection} 
\begin{align*}
 {P}_{\cL,\rho} (s, t) 
 &= \sum_{w\in \cL} \rho(w)  =\sum_{n,d \in \N_0^2} \sum_{w\in \cL_n^d} \rho(w) 
  =\sum_{n,d \in \N_0^2}\#(\cL^d_n) \cdot t^ds^n\\
 &=\sum_{n,d \in \N_0^2} \#[\Mon (A_{n+1})]_d  \cdot t^d s^n =\sum_{n,d \in \N_0^2} \dim_{\K} [A_{n+1})]_d  \cdot t^d s^n\\
 &= s^{-1} \sum_{n \ge 0 }H_{A_{n+1}} (t) \cdot s^{n+1}\\
 &= s^{-1}  {\rm equivH}_{\mathscr{A }} (s, t).
\end{align*}

Notice that $\cL$ can be written as 
\begin{align}
    \label{eq:languageC}
\cL = \{\tau, \alpha_0, \alpha_1 \tau, \ldots, \alpha_c \tau^c \}^* 
\big \{w \in \{\epsilon, \alpha_0,\ldots,\alpha_c \} %\cup \bigcup_{0 \le i \le c} \alpha_i 
\big \} \{\tau \}^*,  
\end{align}
where $\epsilon$ denotes the empty word, 
This shows that $\cL$ is a regular language. % (see also \Cref{exa:easy language}). 
Hence 
${\rm equivH}_{\mathscr{A }} (s, t)$  is a rational function (see \Cref{eq:automaton-HS}). 
\end{proof}

As a regular language,  $\cL$ is recognizable by a finite automaton, which can be used to  compute explicitly the rational form of ${\rm equivH}_{\mathscr{A }} (s, t)$.  

\begin{ex} 
    \label{exa:automaton for filtration}
The regular language in \Cref{eq:languageC} is recognized by a finite automaton, which we describe as a directed labeled graph with nodes/states $(i, j) \in \Z^2$ with $0 \le j \le i \le c$, initial state $(0, 0)$ and accepting states $(i, 0)$ and $(i, i)$ for $i \in \{0,\ldots,c\}$, and edges 
\begin{itemize}
    \item from $(i, i)$ to $(i, i)$ labeled by $\tau$ if $i \in \{0,\ldots,c\}$, 
    \item from $(0, 0)$ to $(i, 0)$ labeled by $\alpha_i$ if $i \in \{0,\ldots,c\}$,
    \item from $(i, j)$ to $(i, j+1)$ labeled by $\tau$ if $0 \le j < i \le c$,
    \item from $(i, i)$ to $(j, 0)$ labeled by $\alpha_j$ if $i, j \in \{0,\ldots,c\}$. 
\end{itemize}
For $c=2$, the automaton is depicted in \Cref{fig:c2automaton}, where $i j$ denotes state $(i, j)$. 
%\vspace*{-1.2cm}
\begin{figure}[h]
    \centering
\tiny
\begin{tikzpicture}[shorten >=1pt,node distance=3cm,auto]
  \tikzstyle{every state}=[fill={rgb:black,1;white,10}]
   \node[state,initial,accepting]  (00)   {$00$};
  \node[state,accepting]  (10) [above right of=00]    {$10$};
  \node[state,accepting]  (11)  [ right of=10]  {$11$};
  \node[state,accepting]  (20) [ below right of=00]   {$20$};
  \node[state]  (21)  [right of=20]  {$21$};
  \node[state,accepting]  (22)  [ right of=21]  {$22$};
 \path[->]
    (00)   edge  [loop above]    node {$\tau,\alpha_0$} (00)
    (00)   edge  []    node {$\alpha_1$} (10)
    (00)   edge  []    node {$\alpha_2$} (20)
    (10)   edge  []    node {$\tau$} (11)
    (20)   edge  []    node {$\tau$} (21)
    (21)   edge  []    node {$\tau$} (22)
    (11)   edge  []    node {$\alpha_0$} (00)
    (11)   edge  [bend right]    node [left]{$\alpha_1$} (10)
    (11)   edge  []    node {$\alpha_2$} (20)
    (22)   edge  []    node {$\alpha_0$} (00)
    (22)   edge  []    node {$\alpha_1$} (10)
    (22)   edge  [bend left]    node {$\alpha_2$} (20)
    (11)   edge  [loop above]    node {$\tau$} (11)
    (22)   edge  [loop above]    node {$\tau$} (22)
    ;
\end{tikzpicture}
\caption{\small\emph{The finite automaton in \Cref{exa:automaton for filtration} for $c=2$.}}
    \label{fig:c2automaton}
\end{figure}
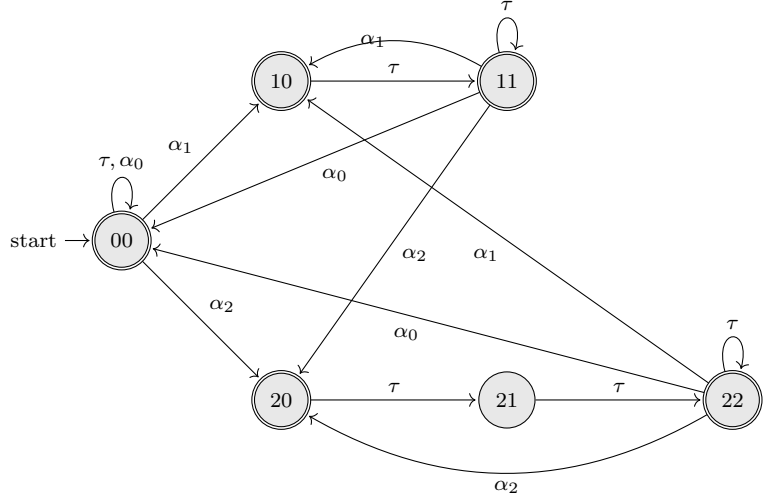

Using \Cref{eq:automaton-HS}, we get for $c=2$, 
\begin{equation*}
P_{\cL, \rho} (t,s)=\begin{bmatrix}
    1\\1\\1 \\1\\0 \\ 1
 \end{bmatrix}^T \cdot 
 \begin{blockarray}{cccccc}
   \textcolor{gray}{00} & \textcolor{gray}{10} & \textcolor{gray}{11}& \textcolor{gray}{20} & \textcolor{gray}{21} & \textcolor{gray}{22} \\
 \begin{block}{[cccccc]}
 1-t-s &0 &-t  &0 & 0 & -t  \\
    -t &1 &-t  &0 & 0 & -t  \\
    0 &-s& 1-s &0 & 0 & 0 \\
    -t&0 &-t  &1 & 0 & -t \\
    0 &0 &0  &-s& 1 & 0 \\
     0 &0 &0  &0& -s & 1-s\\
      \end{block} 
%   \textcolor{gray}{00} & \textcolor{gray}{10} & \textcolor{gray}{11}& \textcolor{gray}{20} & \textcolor{gray}{21} & \textcolor{gray}{22} \\      
\end{blockarray}^{\raisebox{-3ex}{\footnotesize -1}}  \cdot
\begin{bmatrix}
    1\\0\\0 \\0\\0 \\ 0\\0
    \end{bmatrix} \\
    = \dfrac{2t+1}{1-s -ts-ts^2}.
\end{equation*}
Further computations suggest that, for any $c \ge 1$ in \Cref{thm:thm1}, one has $ P_\cL(t,s)= \dfrac{ct+1}{1-s -t\sum\limits_{i=1}^cs_i}$, and so 
\begin{align*}
{\rm equivH}_{\mathscr{A }} (s, t) = \dfrac{s (ct+1)}{1-s -t\sum\limits_{i=1}^cs_i}.
\end{align*}
Code for computing this series for any fixed $c$ can be obtained at  \url{https://sites.google.com/view/aidamaraj/research}.
\end{ex}{}

We conclude this section by showing that Segre products of the above filtrations also have rational equivariant Hilbert series. Consider another filtration $\sB=(B_n)_{n\in \N}$ of the same type as $\cA$, that is,  $
B_n = \K[x_i^2,x_ix_{i+1},\dots,x_ix_{i+c'} \:\mid\: i \in  [n] ] 
$ 
for some integer $c' \ge 0$. 

\begin{thm} 
    \label{thm:Segre product with rat hilb}
The Segre product of the 1-filtrations $\sA$ and $\sB$ has a rational equivariant Hilbert series, that is,  
\[
{\rm equivH}_{\sA \boxtimes \sB} (s_1,s_2, t) = \sum_{m, n \in \N} H_{A_m\boxtimes B_n} (t) \cdot s_1^m s_2^n = 
\sum_{m, n \in \N} \sum_{d \ge 0} \dim_{\K} [A_m \boxtimes B_n]_d \cdot t^d s_1^m s_2^n 
\]
is a rational function. 
\end{thm}

\begin{proof} 
By \Cref{lem:bijection}, the filtrations $\sA$ and $\sB$ are represented by languages that are regular (see Equation \eqref{eq:languageC}).  Hence we conclude by \Cref{prop:relate Segre products}. 
\end{proof}

%%%%%%%%%%%%%%%%%%%%%%%%%%%%%%%%%

\section{An Infinitely Generated Toric Presentation Ideal}
    \label{sec:infinitely gen}

In \Cref{sec:1} we considered the shift invariant subalgebras of $\K[X_{\N}] = \K[x_i \; \mid \; i \in \N]$: 
\[
\K[x_i^2, x_i x_{i+1} \; \mid \; i \in \N] \quad  \text{ and } \quad  \K[x_i^2, x_i x_{i+1}, x_i x_{i+2}  \; \mid \; i \in \N]. 
\] 
We showed that both algebras are the limits of families of algebras with a rational Hilbert series. As a consequence of a later more general result (see \Cref{rem:presentation ideal}),  both algebras have presentation ideals that are finitely generated up to shift. These ideals are the kernels of the shift invariant homomorphisms
\[
\K[X_{[2]\times \N}] = \K[x_{1, j}, x_{2, j} \; \mid \; i \in \N] \to \K[X_{\N}], \ x_{1, j} \mapsto x_j^2, x_{2, j} \mapsto x_j x_{j+1}
\]
and 
\[
\K[X_{[3]\times \N}] = \K[x_{1, j}, x_{2, j}, x_{3, j}  \; \mid \; i \in \N] \to \K[X_{\N}], \ x_{1, j} \mapsto x_j^2, x_{2, j} \mapsto x_j x_{j+1},  x_{3, j} \mapsto x_j x_{j+2}, 
\]
respectively. In this section, we mainly consider a seemingly similar monomial subalgebra, namely  
\[
A = \K[x_i x_{i+1}, x_i x_{i+2}  \; \mid \; i \in \N] \subset \K[X_{\N}]. 
\]
It is the image of the shift invariant monomial homomorphism
\begin{align}\label{eq:mainex}
    \varphi \colon \K[X_{[2]\times \N}]  \longrightarrow \K[X_{\N}], \ x_{1,i} \mapsto x_i x_{i+1}, \ x_{2,i} \mapsto x_{i}x_{i+2}.
\end{align} 
However, in this case we show that the presentation ideal $I =  \ker(\varphi)$ is \emph{not} finitely generated up to shift. In fact, we explicitly describe a minimal generating set of $I$ up to shift  that consists of binomials of arbitrarily large degree. 
\smallskip 

To establish this result, it is convenient to change notation. We will write $x_{i,i+1}$ instead of $x_{1, i}$ and $x_{i,i+2}$ instead of $x_{2, i}$. Thus, the above map $\varphi$ becomes the surjective ring homomorphism
\[
\varphi \colon R\coloneqq K[x_{i,i+1},x_{i,i+2} \;  \mid \;  i\in \N] \to A, \text{ with } 
x_{i,i+1} \mapsto x_i x_{i+1}, \ x_{i,i+2} \mapsto x_{i}x_{i+2}.
\]
Denote by $E=\{[i,i+1], \ [i,i+2] \; \mid \; i\in \N\}$ the set of index pairs of variables in $R$. Any binomial in $R$ can be written as $\bx^\bu-\bx^\bv  = \prod_{[i,j]\in E}x_{i,j}^{u_{i,j}}-\prod_{[i,j]\in E}x_{i,j}^{v_{i,j}}$ with non-negative integers $u_{i, j}$ and $v_{i, j}$. Abusing notation slightly, we write $k \in [i, j]$ if $k = i$ or $k = j$. It follows that 
\begin{align}
    \label{eq:kernel}
    \bx^\bu-\bx^\bv \in I = \ker \varphi \text{ if and only if, for any $k \in N$, one has } \sum\limits_{k \in [i, j]} u_{i,j}=\sum\limits_{k \in [i, j]} v_{i,j}.
\end{align} 
    
We now describe binomials in $R$ by weighted graphs. To this end 
let $G$ be the graph  with vertex set $\N$ and edge set $E$. Consider any  \emph{weight function}  $w\colon E\rightarrow \Z, \ [i,j]\xrightarrow[]{} w([i,j])$ on the edges of $G$.  It assigns weights to the vertices of $G$ by using adjacent edges, that is,   $w(n)=\sum\limits_{n \in [i,j]} w([i,j])$ for any vertex $n$ of $G$.  Explicitly, this gives  
\begin{align}
    w(1)&=w([1,2])+w([1,3]), \label{eq:weights1}\\
    w(2)&=w([1,2])+w([2,3])+w([2,4]), \text{ and } \label{eq:weights2} \\
   w(i)&=w([i-2,i])+w([i-1,i])+w([i,i+1])+w([i,i+2]) \text{ if } i \ge 3. \label{eq:weights3}
\end{align}

We encode any binomial $g=\bx^\bu-\bx^\bv$ in $R$ by a weighted graph $G_g$ with vertex set $\N$ whose edges are the edges of $G$ in the support of $\bu$ or $\bv$ 
We weigh the edges of $G_g$ 
using the weight function $w_g \colon E \to \Z$, defined by $w_g([i,j])=u_{i,j}-v_{i,j}$. Thus, the edges of $G_g$ are precisely the edges $[i, j]$ of $G$ with non-zero weight $w_g([i,j])$  (see \Cref{fig:g4} for pictures of some graphs).  
Note that any two binomials $g$ and $x_i g$ are encoded by the same weighted subgraph of $G$. However, in order to describe a generating set of $I = \ker \varphi$ it is enough to consider binomials $\bx^\bu-\bx^\bv$ that are not a multiple of any variable, that is, the supports of $\bu$ and $\bv$ are disjoint. Any binomial with this property is uniquely identified by its weighted graph. In any case, this graph can be used to decide whether the binomial belongs to $I$. 

\begin{prop}
    \label{lm:bijection}
A binomial $g=\bx^\bu-\bx^\bv$ is in $I$ if and only if the weight function $w_g$ satisfies $w_g(n)=0$ for every vertex $n\in \N$. 
\end{prop}

\begin{proof}
This is an immediate consequence of Condition \eqref{eq:kernel}. 
\end{proof}

Since the ideal $I$ is shift invariant, \Cref{lm:bijection} and \Cref{eq:weights1} give the following observation. 

\begin{cor}
     \label{cor:smallest_vertex}
Consider any binomial $g=\bx^\bu-\bx^\bv$ in $I$ that is not a product of a variable,  and denote by $c$ the smallest vertex adjacent to an edge of its weighted graph $G_g$. Then $w_g([c,c+1])=-w_g([c,c+2]) \neq 0$. 
\end{cor}

Now, we define recursively a set of binomials which, we will show, minimally  generates~$I$. For a finite sequence $\bs = (s_1,\ldots,s_n) \subset \{1,2\}^n$ with $n\in \N_0$ , we set $\bs 1 = (s_1,\ldots,s_n, 1)$
and $\bs 2 = (s_1,\ldots,s_n, 2)$. We define $\{1,2\}^0 = \emptyset$. 
Observe that the partial order on $\N^2$ defined by $(i, j) \le (k, l)$ if $i \le k$ and $j \le l$ induces a total order on the edge set $E$ of $G$. 
Recall that any binomial $g=\bx^\bu-\bx^\bv \in R$ that is not a multiple of a variable and its graph $G_g$ are defined by a weight function $E \to \Z$ with finite support. 

\begin{defn}
   \label{def:generating set} 
Define recursively a set of binomials  $\G' = \{ g^\bs \; \mid \; \bs  \subset \{1,2\}^n, \ n \in \N_0 \}$     
by 
\[
g^{\emptyset}=x_{1, 2}x^2_{3, 5}x_{6,7}-x_{1,  3}x_{2, 3}x_{5,6}x_{5,7}
\]
and, if $g^{\bs}$ is defined and $k$ is the largest index of a vertex adjacent to an edge of $G_{g^{\bs}}$, 
\begin{align}
    w_{g^{\bs 1}}([i,j])=\begin{cases}
     w_{g^{\bs}}([i,j]) & \text{ if } (i,j) \le (k-4,k-2),\\
     2 w_{g^{\bs}}([i,j]) & \text{ if } (i,j) = (k-2,k),\\
    - w_{g^{\bs}}([i-2,j-2]) & \text{ if } (i,j) \ge (k,k+1),\\
    0 & \text{ otherwise},  
    \end{cases}
\end{align}
 and 
\begin{align}
w_{g^{\bs 2}}([i,j])=\begin{cases}
w_{g^s}([i,j]) & \text{ if } (i,j) \le (k-4,k-2),\\
2 w_{g^{\bs}}([i,j]) & \text{ if } (i,j) = (k-2,k-1),\\
 - 2 w_{g^s}([k-2,k-1]) & \text{ if } (i,j)=(k-1,k+1),\\
w_{g^s}([i-3,j-3]) & \text{ if } (i,j) \ge (k+1,k+2), \\
0 & \text{ otherwise}. 
    \end{cases}
\end{align}

Finally, set 
\[
\G = \{ g_2\} \cup \G', \quad \text{ where } \; g_2=x_{1, 2}x_{3, 4}-x_{1, 3}x_{2, 4}. 
\]
\end{defn}

\begin{ex}
We illustrate the passage from $g^{\emptyset}$ to $g^{(1)}$ and $g^{(2)}$ in  \Cref{fig:g4}. 
\begin{figure}[ht]
    \centering
    \includegraphics[scale=0.55]{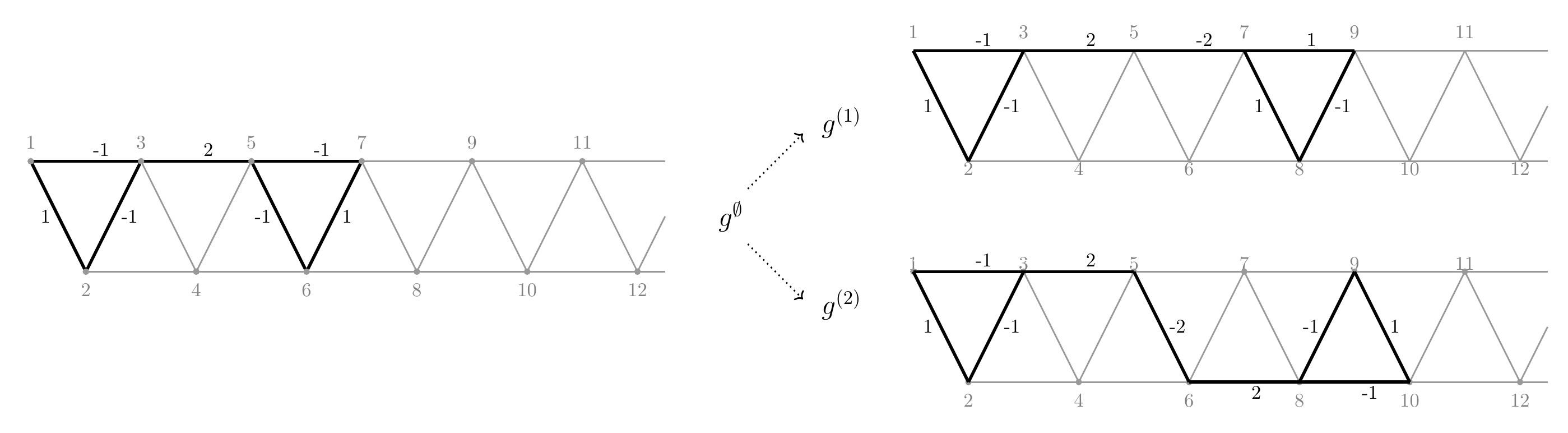}
    \caption{Visualizations of $g^\emptyset$, $g^{(1)}$ and $g^{(2)}$. }
    \label{fig:g4}
\end{figure}
\end{ex}

The binomials in $\G'$ and their weighted graphs have the following properties. 

\begin{prop}
    \label{prop:properties of generating binomials} 
For any $g \in \G'$, one has: 
\begin{itemize}

\item[(a)] Let $k$ be the largest index of a vertex adjacent to $G_g$. Then  one has $k \ge 7$ and the support of $G_g$ consists of two triangles $[1, 2], [1, 3], [2, 3]$ and $[k-2, k-1], [k-2, k], [k-1, k]$ with $w_g([1, 3]) = w_g ([2, 3] ) = - w_g ([1, 2]) = - 1$ and $w_g ([k-2, k-1]) = w_g ([k-2, k]) = - w_g ([k-1, k]) \in \{1, -1\}$ that are connected by a path $[3,5] = [i_1, i_2]\ldots,[i_{l-1}, i_l] = [k-4, k-2]$ with $(i_1, i_2) < \cdots < (i_{l-1}, i_l)$ whose edges have weights alternating between 2 and $-2$ and no consecutive edges of the form $[j-1,j], [j, j+1]$.   

\item[(b)] $\deg g =4+\sum_{i=1}^n s_i$ if $g = g^{\bs}$ with $\bs = (s_1,\ldots,s_n)$. 

\item[(c)] $g \in I = \ker \varphi$. 
\end{itemize}

\end{prop}

\begin{proof}
Claim (a) follows by analyzing the recursive definition of $\G'$.
For (b), note that $\deg g^{\bs 1} = 1 + \deg g^{\bs}$ and $\deg g^{\bs 2} = 2 + \deg g^{\bs}$. 
Claim (c) is a consequence of (a) and Condition \eqref{eq:kernel}. 
\end{proof}

Recall that the sequence of Fibonacci numbers $(F_n)_{n \ge 0}$ is defined by $F_0 = 0, F_1 = 1$ and $F_{n+2} = F_{n+1} + F_n$ if $n \in \N_0$. 

\begin{cor}\label{prop:fibonacci}
The number of binomials of degree $d\geq 3$ is $\G'$ is equal to the Fibonacci number $F_{d-3}$.
\end{cor}

\begin{proof}
The number of sequences $\{1, 2\}^n$ with $n \in \N_0$ whose entries sum up to $k \ge 0$ is equal to $F_{k+1}$. Hence we conclude by \Cref{prop:properties of generating binomials}(b). 
\end{proof}

The following observation will be useful. 

\begin{lm}
   \label{lem:Markov move}
Consider a prime ideal $I$ that is generated by binomials and a binomial $h = \bx^{\ba} - \bx^{\bb} \in I$. 
Then one has: 
\begin{itemize}

\item[(a)] If there is some binomial $g=\bx^\bu-\bx^\bv \in I$ with $\ba \ge \bu$ and $\bb \ge \bv - \bv' \ge 0$ for some $\bv' \ge 0$, then 
\[
\bx^{\ba - \bu + \bv'} - \bx^{\bb - (\bv - \bv')} 
\]
is in $I$. 

In particular, if $h$ and $g$ are homogeneous and $\bv' < \bv$, then this binomial has smaller degree than $h$. 

\item[(b)] If $B$ is a generating set of $I$ consisting of binomials, then there is some $g=\bx^\bu-\bx^\bv \in B$ such that 
$\ba +\bu - \bv \ge 0$ or $\ba - \bu+\bv \ge 0$. 
\end{itemize}
\end{lm}

\begin{proof} 
(a) 
Since $h$ and $g$ are in $I$, so is 
\begin{align*}
h - \bx^{\ba - \bu} \cdot g & = \bx^{\ba} - \bx^{\bb} - \bx^{\ba -\bu} \big ( \bx^\bu-\bx^\bv  \big ) \\
&  = \bx^{\bv - \bv'} \big ( \bx^{\ba - \bu + \bv'} - \bx^{\bb - (\bv - \bv')} \big ). 
\end{align*}
Using that $I$ is a prime ideal, the claim follows. 

(b) Since $B$ consists of binomials and generates $I$, one can write $h \in I$ as 
\[
\bx^{\ba} - \bx^{\bb} = \sum_{k = 1}^l \varepsilon_k \bx^{\bm_k} (\bx^{\bu_k} - \bx^{\bv_k}) 
\]
with binomials $\bx^{\bu_k} - \bx^{\bv_k} \in B$, monomials $\bx^{\bm_k}$ and $\varepsilon_k \in \{1, -1\}$. Possibly after re-indexing, it follows that $\ba = \bm_1 + \bu_1$ or $\ba = \bm_1 + \bv_1$. In the former case, we get
\[
\bx^{\ba} - \bx^{\bb} - \bx^{\bm_1} (\bx^{\bu_1} - \bx^{\bv_1}) =   \bx^{\bm_1 + \bv_1} - \bx^{\bb} =  \bx^{\ba - \bu_1 + \bv_1} - \bx^{\bb}, 
\]
which shows $\ba - \bu_1 + \bv_1 \ge 0$. In the latter case, we obtain similarly $\ba + \bu_1 - \bv_1 \ge 0$. 
\end{proof}

We now establish the main result of this section. 

\begin{thm} 
    \label{thm:min gen set}
The set $\G=\{g_2,g_4^\bs \mid \bs\in \{1,2\}^n, n\in \N_0\}$ generates $I$ minimally up to shift. 
\end{thm}

\begin{proof}
Since $I$ is the kernel of a monomial map, it has a minimal generating set consisting of binomials. 

First, we show that $\G$ generates $I$ up to shift. One easily checks that $g_2$ is in $I$. Combined with \Cref{prop:properties of generating binomials}, this gives $\G \subset I$. It remains to show $I \subseteq \langle \sh(\G) \rangle$, where $\langle \sh(\G) \rangle$ is the ideal generated by shifts of the polynomials in $\G$. 

As $I$ is shift invariant it suffices to prove that any binomial $h \in I$ with the property that vertex $1$ is adjacent to an edge of $G_h$ is in the ideal generated by $\G$.  Let $h = \bx^{\bu}-\bx^{\bv} $ be a binomial of least degree that is not in in the ideal generated by $\G$. Thus,  $h$ is not a multiple of any variable.  By considering several cases, we will argue that $h$ can be used to produce a binomial in $I \setminus \langle \sh(\G) \rangle$ whose degree is less than $d = \deg h$. This contradiction will prove that $I = \langle \sh(\G) \rangle$. 

For simplicity let us write $w$ for the weight $w_h$ defined by $h$. By Equation \eqref{eq:weights1} and \Cref{lm:bijection}, we must have that $w([1,2])=-w([1,3]) \neq 0$. Replacing possibly $h$ by $-h$, we may assume $w([1,2])>0$. Thus, Equation \eqref{eq:weights2} and \Cref{lm:bijection} imply  $w([2,3])+w([2,4])<0$.   If $w([2,4])<0$, then we use \Cref{lem:Markov move} with $g = g_2$ to obtain a binomial that is not in $\langle \sh(\G) \rangle$, as desired. 

Notice that this last argument gives more generally that the weighted graph $G_h$ cannot have, for any $k \in \N$,  three out of four edges $[k, k+1], [k, k+2], [k+1, k+3], [k+2, k+3]$ with the property that $[k, k+1]$ and  
$[k+2, k+3]$ have the same nonzero sign opposite to the signs of $[k, k+2]$ and  $[k+1, k+3]$. Otherwise, \Cref{lem:Markov move} applies as above using a suitable shift of $g_2$. We will use this observation in the remainder of the proof. 

As a first consequence, we conclude $w([2,4]) \ge 0$ and $w(3,4) \le 0$, which implies $w([2,3])<0$ and,   by using $w([3,4])+w([3,5])=-(w([1,3])+w([2,3]))\geq 2$, also  $w([3,5]) \ge 2$. The latter gives $w([4,5]) \ge 0$ because otherwise \Cref{lem:Markov move} with $-\sh_1(g_2)=x_{2,4}x_{3,5}-x_{2,3}x_{4,5}$ is applicable. 

% start of induction
Let us say that  $h \in I$ has \emph{Property $(+)$} if, for some $k \ge 5$,  
\begin{itemize}

\item[(i)] its graph $G_h$ contains the triangle $[1, 2], [1, 3], [2, 3]$ and a path $[3,5] = [i_1, i_2]\ldots,$ 
$ [i_{l-1}, i_l] = [k-2, k]$ with $(i_1, i_2) < \cdots < (i_{l-1}, i_l)$ whose edges have alternating weights with absolute value at least two and no consecutive edges of the form $[j-1,j], [j, j+1]$; \; and 

\item[(ii)]  $w([k-1,k])$ and $-w([k-2, k])$ do not have the same (non-zero) sign. 
\end{itemize}

Notice that we showed above that $h$ satisfies Property $(+)$ with $k=5$. Our goal is to show the following \emph{Claim}:  If $h$ satisfies Property $(+)$ for some $k \ge 5$, then either one can replace $h$ by a binomial of smaller degree by using  \Cref{lem:Markov move}, as desired, or $h$ satisfies Property $(+)$ with some $k' > k$. 

As $h$ involves only finitely many variables, the Claim implies that eventually we must arrive in the first situation, which gives the desired contradiction to the minimality of $\deg h$ among the polynomials in $I \setminus \langle \sh(\G) \rangle$. 

It remains to establish the above Claim. Assume $h$ has Property $(+)$ for some $k \ge 5$. Possibly replacing $h$ by $-h$, we may assume $w([k-2, k]) \ge 2$ and $w([k-1,k]) \ge 0$.  If $w([k,k+1])<0$ and $w([k,k+2])<0$, then one can use \Cref{lem:Markov move} with the element of $\G'$ that contains the triangle $[k,k+1], [k, k+2], [k+1, k+2]$ and the path from vertex $5$ to vertex $k$ as described in Property $(+)$. This shows  the Claim in this case. Otherwise, the conditions imposed on vertex $k$ by  \Cref{lm:bijection} show that one of the following two conditions must be satisfied: (a)  $w([k,k+1])\leq -2$ and $w([k,k+2])\geq 0$, or  (b) $w([k,k+1])\ge 0$ and $w([k,k+2])\le -2$. 
We consider these separately. 

\emph{Case (a)}: Note that  we must have $w([k-1, k+1]) \le 0$ because otherwise one can reduce $h$ using a shift $g_2$ to a binomial of smaller degree. Similarly, we get $w([k+1,k+2]) \le 0$ as otherwise one one can reduce $h$ using the element of $\G'$ that contains the triangle $[k, k+1], [k, k+2], [k+1, k+2]$ and the path from vertex $5$ to vertex $k$ as described in Property $(+)$.   
Hence, the condition on vertex $k+1$ gives $w([k+1,k+3]) = - w([k-1,k+1]) - w([k+1,k+2]) - w([k,k+1]) \ge 2$. It follows now that $w([k+2, k+3]) \ge 0$ as otherwise one can reduce $h$ using a shift $g_2$. The last two inequalities show that $h$ satisfies Property $(+)$ for $k+3$. 

\emph{Case (b)}: Observe that  we must have $w([k+1, k+2]) \le 0$ because otherwise one can reduce $h$ using the element of $\G'$ that contains the triangle $[k, k+1], [k, k+2], [k+1, k+2]$ and the path from vertex $5$ to vertex $k$ as described in Property $(+)$. Since we know by the assumption for (b) that $w([k,k+2])\le -2$, we see that $h$ satisfies Property $(+)$ for $k+2$. 

This completes the proof of the Claim and establishes that $\G$ generates $I$ up to shift. 
\smallskip

Second, we establish that $\G$ is a  minimal generating set of $I$ up to shift. To this end we prove that for any $h \in \G$, there is no $\tilde{g} \in \G \setminus \{h \}$ such that $h$ and any shift $g$ of $\tilde{g}$ satisfy one of the two conditions in \Cref{lem:Markov move}(b). This follows once we have shown that, for any such $h$ and $g$,  there are always  $[i,j], [k,l] \in E$ such that $w_g([i,j])= w_g([k,l])=0$ but $w_h ([i,j])>0$ and $w_h([k,l])<0$. 

Suppose $g$ is properly shifted, that is, vertex $1$ is not adjacent to any edge of  the graph $G_g$, and so $w_g([1,2])= w_g([1,3])=0$. 
Since $- w_h([1, 3]) =  w_h ([1, 2]) = 1$ by \Cref{prop:properties of generating binomials}(a), we are done in this case. 

Thus, we may assume $g = \tilde{g}  \in \G \setminus \{h \}$. Let $m$ be the maximal vertex that is adjacent to $G_h$. 
If $m$ is not the maximal vertex that is adjacent to $G_g$, we use an argument similar to the one in the previous graph. 
Indeed, in this case we may assume that the largest vertex that is adjancent $G_g$ is less than $m$, and so 
$w_g([m-2,m])= w_g([m-1,m])=0$. However, by the choice of $m$, \Cref{lm:bijection} at vertex $m$ gives $w_h([m-2,m])= - w_h([m-1,m) \neq 0$, as desired. 

Hence, we are left to consider $h \neq g$ in $\G$ such that $m$ is also the maximal vertex that is adjacent to $G_g$. Thus, both $G_h$ and $G_g$ contain a path from vertex $5$ to vertex $m-4$ as described in \Cref{prop:properties of generating binomials}(a). It the two paths have the same support then we get $h = g$. Thus, there is a maximum vertex $c \ge 5$ such that the paths from $5$ to $c$ in $G_h$ and $G_g$ are the same. We may assume that the path in $G_h$ continues with the edge $[c, c+1]$ and the next edge in $G_g$ is $[c, c+2]$. By \Cref{prop:properties of generating binomials}(a), it follows that vertex $c+1$ is not adjacent to $G_g$, which gives on the one hand
\[
w_g ([c, c+1]) = w_g ([c+1, c+2]) = w_g ([c+1, c+3]) = 0. 
\]
On the other hand, \Cref{prop:properties of generating binomials}(a) yields $w_h ([ c, c+1]) = - w_h ([ c+1, c+2]) \neq 0$ or $w_h ([ c, c+1]) = - w_h ([ c+1, c+3]) \neq 0$.  This completes the argument. 
\end{proof}

\begin{cor}
    \label{cor:ideal not fg}
The ideal $I$ is not finitely generated up to shift. 
\end{cor} 

\begin{proof}
The degrees of the binomials in the minimal generating set $\G$ are not bounded. 
\end{proof}

In contrast, the monomial algebras considered in \Cref{sec:1} have presentation ideals that are finitely generated up to shifting. 

\begin{prop}
     \label{rem:presentation ideal}
Fix any integer $c \ge 1$ and consider the $\K$-algebra  $A=\K[x_i x_j \; \mid \; i, j  \in \N, \  i \le j \le i+c ]$ and the surjective, shift-equivariant homomorphism $\varphi \colon \K[X_{[c+1]\times \N}]  \longrightarrow A$, defined by $x_{i, j} \mapsto x_j x_{i + j -1}$.  The presentation ideal $I = \ker \varphi$ of $A$ is finitely generated up to shifting by quadrics. 
\end{prop}

\begin{proof} 
As above, it is convenient to change notation. We write $x_{j,i+j}$ instead of $x_{i+1, j}$. Thus, $\varphi$ becomes the homomorphism
\[
\psi \colon \K[x_{i, j} \; \mid \; i, j \in \N, \ 0 \le j-i \le c] \to A \; \text{ with } x_{i, j} \mapsto x_i x_j. 
\]
Note that shifting acts on the domain by $sh_k (x_{i, j}) = x_{i+k, j+k}$. 
We claim that $J = \ker \psi$ is generated by the set $\G'$ of quadrics  
\begin{align*}
x_{i,j}x_{\min\{k,\ell\}, \max\{k,\ell\}}-x_{i,k}x_{\min\{j,\ell\}, \max\{j,\ell\}}
%
%x_{i,j}x_{k,\ell}-x_{i,k}x_{j,\ell} \quad \text{ with } 0\leq j-i,\ell-k,k-i,\ell-j\leq c, \\
%x_{i,j}x_{k,\ell}-x_{i,k}x_{\ell,j} \quad \text{ with }  0\leq j-i,k-\ell,k-i,\ell-j\leq c, \\
%x_{i,j}x_{k,\ell}-x_{i,\ell}x_{j,k} \quad  \text{ with }  0\leq j-i,k-\ell,\ell-i,k-j\leq c, \\
% x_{i,j}x_{k,\ell}-x_{i,\ell}x_{k,j} \quad \text{ with }  0\leq j-i,k-\ell,\ell-i,k-j\leq c.
\end{align*}
with $i \le j < k \le i+c$, $i < \ell$, $|k - \ell| \le c$ and $|j - \ell| \le c$. 
This will prove the claim because any quadric in $\G'$ can be obtained by shifting one of the above quadrics with $i =1$, and there are only finitely many quadrics in $\G'$ with $i=1$. 

To establish the claim note that $\G'$ is in $J$. It remains to show that $\G'$ generates $J$. 
Let $h=\bx^{\bu}-\bx^{\bv}$ be a binomial in $J$ of least degree that is not in the ideal generated by $\G'$. Thus $\bu$ and $\bv$ have disjoint support. Order $\N_0^2$ lexicographically, that is, $(i, j) < (k, \ell)$ if $i < k$ or $i =k$ and $j < \ell$. 
Let $(i,j)$ be the smallest index of a variable appearing in the support of $h$. We may assume that it appears in $\bx^{\bu}$,  and so $u_{i,j}\neq 0$. Since $\varphi(\bx^{\bu})=\varphi(\bx^{\bv})$, there is some $k>j$  such that $v_{i,k}\neq 0$ and some $\ell> i$ with  $0\leq |\ell-j|\leq c$ such that $v_{\min \{ \ell,j \},\max \{ \ell,j \}}\neq 0$.  There are two cases. One has (i) $i < \ell \le j$ and $\ell \ge j$ or (ii)  $i < \ell$ and $\ell \ge j$. In both cases one checks that the quadric 
$q = x_{i,j}x_{\min\{k,\ell\}, \max\{k,\ell\}}-x_{i,k}x_{\min\{j,\ell\}, \max\{j,\ell\}}$ is in the domain of $\psi$. In fact, $q$ is in $J$. 
If $h$ has degree two then it must be equal to $q$ up to sign. If $\deg h \ge 3$ we use \Cref{lem:Markov move} with $g = q$ to obtain a binomial that is not in $\langle \G' \rangle$ and whose degree is less than $\deg h$. This contradicts the choice of $h$ and completes the argument. 
%
%If $\ell<k$, then $0\leq |k-\ell|\leq |k-i|\leq c$.  If  $\ell>k$, then $0\leq |k-\ell|\leq |j-\ell|\leq c$. In both scenarios, $0\leq |k-\ell|\leq c$  and \(g=x_{i,j}x_{\max\{k,\ell\},\min\{k,\ell\}}-x_{i,k}x_{\min\{j,\ell\}, \max\{j,\ell\}}\) is one of the bimomials  in $\ker(\varphi)$ described earlier. 
%\Cref{lem:Markov move}(a) reduces the problem to a binomial of lesser degree.
\end{proof}

%%%%%%%%%%%%%%%%%%%%%%%%%%%%%%%%%%%%%%%%%%%%%%%%

\section{A Rational Hilbert Series of an Infinitely Generated Ideal} 
\label{sec:Hilb non-fg ideal}

Using the results of \Cref{sec:infinitely gen}, we discuss a filtration of monomial algebras $\sA=(A_n)_{n\in \N}$ whose colimit
$A = {\displaystyle \lim_{\longrightarrow}}\, A_n$ is finitely generated up to shift. These algebras have toric presentation ideals, denoted $I_n$ and $I = {\displaystyle \lim_{\longrightarrow}}\, I_n$. The ideal $I$ is shift invariant, but not finitely generated up to shift by \Cref{cor:ideal not fg}. Nevertheless, we show that its equivariant Hilbert series 
$\equivH_I (s,t)=\sum\limits_{n\geq 1,d\geq 0} \dim_{\K} [I_n]_dt^ds^n$ is rational (see \Cref{cor:Hilb of non-fg ideal}). To the best knowledge of the authors this is the first example of this kind. Previous rationality results for equivariant Hilbert series were established for finitely generated modules (see, e.g., \cite{N}).  

We continue to use the notation introduced at the beginning of \Cref{sec:infinitely gen}, that is, we consider the subalgebra 
$A=\K[x_i x_{i+1}, x_i x_{i+2}  \mid i \in \N]$ of $\K[X_{\N}]$. It is the image of the map
%see Equation \eqref{eq:mainex}
\begin{align*}
    \varphi \colon \K[X_{[2]\times \N}]  \longrightarrow \K[X_{\N}], \ x_{1,i} \mapsto x_i x_{i+1}, \ x_{2,i} \mapsto x_{i}x_{i+2}.
\end{align*} 
Its presentation ideal $I = \ker \varphi$ is shift invariant, but not finitely generated up to shift (see \Cref{cor:ideal not fg}). Consider a filtration $\sA=(A_n)_{n\in \N}$ of subalgebras, where $A_n = \K [x_i x_{i+1}, x_i x_{i+2} \mid i \in [n]]$. Define their presentation ideals as $I_n = \ker \varphi_n$ with 
\begin{align*}
    \varphi_n \colon K[X_{[2]\times [n]}]  \longrightarrow \K[X_{n+2}], \ x_{1,i} \mapsto x_i x_{i+1}, \ x_{2,i} \mapsto x_{i}x_{i+2}, 
\end{align*} 
where $K[X_{[2]\times [n]}]  = \K[x_{i, j} \mid i \in [2], j \in [n]]$ and $\K[X_{n+2}] = \K[x_i \mid i \in [n+2]]$. Notice that 
$\K[X_{[2]\times \N}]/I \cong  A = {\displaystyle \lim_{\longrightarrow}}\, A_n \cong {\displaystyle \lim_{\longrightarrow}}\, \K[X_{[2]\times [n]}]/I_n$ and $I = {\displaystyle \lim_{\longrightarrow}}\, I_n$. One can also recover $I_n$ from $I$.

\begin{prop}
\begin{itemize}
\item[(a)] For every $n \in \N$, one has $I_n = I \cap K[X_{[2]\times [n]}]$. 

\item[(b)] The maximum degree of a minimal generator of $I_n$ approaches infinity as $n$ approaches infinity. In fact, if $n \ge 6$ this maximum degree is  $\frac{2n}{3}$ if $n \equiv 0\!\!\! \mod 3$, $\frac{2n+1}{3}$ if $n \equiv 1\!\!\! \mod 3$
and $\frac{2n-1}{3}$ if $n \equiv 2\!\!\! \mod 3$. 
\end{itemize}
\end{prop}

\begin{proof}
(a) Clearly, the extension ideal of $I_n$ in $\K[X_{[2]\times \N}]$ is contained in $I$. Conversely, any polynomial $f \in I$ involves only finitely many variables. If it is contained in $K[X_{[2]\times [n]}]$, then it belongs to $I_n$. 

(b) Similarly, it follows that $I_n$ is minimally generated up to shift by $\G \cap \K[X_{[2]\times [n]}]$, where $\G$ is the minimal generating set of $I$ given in \Cref{thm:min gen set}.  Using also \Cref{prop:properties of generating binomials}(b), the claim follows. 
\end{proof}

Next, we use a suitable formal language to show that the equivariant Hilbert series 
$\equivH_\sA(s,t)=\sum\limits_{n\geq 1,d\geq 0}\dim_{\K}[A_n]_dt^ds^n$ is rational. 

\begin{prop} 
   \label{prop:normal form} 
Every monomial in A has a unique string presentation of the form 
\begin{align}
    \label{eq:string_presentation}
(x_{i_1}x_{i_1+j_1})(x_{i_2}x_{i_2+j_2})\dots (x_{i_d}x_{i_d+j_d})
\end{align}
satisfying 
\begin{enumerate}[(a)]
    \item $i_1\leq \dots i_d$,
    \item $j_1,\dots,j_d\in \{1,2\}$,
    \item if $i_k=i_{k+1}$ then $j_k\leq i_{k+1}$, and 
    \item if $i_{k+1}=i_{k}+1$, then $(j_k,j_{k+1})\neq (2,2)$. 
\end{enumerate}
\end{prop}

\begin{proof}
By definition of $A$ and after possibly changing the order to $(x_x x_{i+2}) (x_i x_{i+1})$, any monomial in $A$ can be written as  $(x_{i_1}x_{i_1+j_1})(x_{i_2}x_{i_2+j_2}) \cdots (x_{i_d}x_{i_d+j_d})$ such that, when read as a string, it satisfies Conditions (a) - (c). Using the order $<$ on $\Z^2$ defined by $(i, j) < (k, l)$ if $i < k$ or if $i = k$ and $j < l$, this means that any monomial in $A$ can be written as 
\begin{align} 
   \label{eq:initial rewriting} 
\bm = (x_{i_1} x_{i_1 + j_1})^{e_1} \cdots  (x_{i_d} x_{i_d + j_d})^{e_d}
\end{align} 
with integers $e_k \in \N$, $j_k \in \{1, 2\}$ and $(i_k, j_k) < (i_{k+1}, j_{k+1})$ whenever $i \le k < d$. We will show that $\bm$ can be rewritten such that it also satisfies Conditions (a) - (d). 

Assume the above description of $\bm$ violates Condition (d) when read as a string. Choose $k$ minimal such that $i_{k+1} = 1 + i_k$ and $j_k = j_{k+1} = 2$. Thus, we can rewrite $\bm$ as
\begin{align} 
       \label{eq:rewritten monomial}
\bm =  \cdots (x_{i_{k-1}} x_{i_{k-1} + j_{k-1}})^{e_{k-1}} \widetilde{\bm} (x_{i_{k}+2} x_{i_{k+2} + j_{k+2}})^{e_{k+2}}  \cdots 
\end{align}
with 
\begin{align*}
\widetilde{\bm}  = 
\begin{cases}
(x_{i_k} x_{i_{k} + 1})^{e_{k+1}}  (x_{i_k} x_{i_{k} + 2})^{e_k -e_{k+1}} (x_{i_k +2} x_{i_{k} + 3})^{e_{k+1}}  & \text{if } e_k \ge e_{k+1}, \\
(x_{i_k} x_{i_{k} + 1})^{e_{k}}  (x_{i_k} x_{i_{k} + 3})^{e_{k+1} -e_{k}} (x_{i_k +2} x_{i_{k} + 3})^{e_{k}} & \text{if } e_k \le e_{k+1}. 
\end{cases}
\end{align*}
Observe that, as a monomial, $\widetilde{\bm} = (x_{i_k} x_{i_{k} + 2})^{e_{k}}  (x_{i_k +1} x_{i_{k} + 3})^{e_{k+1}}$. Consider now the left-most part of $\bm$, that is, 
\[
\bm'  = (x_{i_1} x_{i_1 + j_1})^{e_1} \cdots  (x_{i_{k+1}} x_{i_{k+1} + j_{k+1}})^{e_{k+1}}. 
\]
Since the presentation of $\bm$ in  Equation \eqref{eq:initial rewriting} satisfies  in particular that 
$(i_{k-1}, j_{k-1}) < (i_k, j_k) = (i_k, 2)$, it follows $(i_{k-1}, j_{k-1}) \le (i_k, 1)$.  If this inequality is strict then $\bm$ satisfies Conditions (a) - (d). Otherwise, we rewrite $\bm'$ as the string 
\begin{align*}
\bm' = & \cdots (x_{i_{k-2}} x_{i_{k-2} + j_{k-2}})^{e_{k-2}} \cdot \\
& \begin{cases}
(x_{i_k} x_{i_{k} + 1})^{e_{k-1} + e_{k+1}}  (x_{i_k} x_{i_{k} + 2})^{e_k -e_{k+1}} (x_{i_k +2} x_{i_{k} + 3})^{e_{k+1}}  & \text{if } e_k \ge e_{k+1}, \\
(x_{i_k} x_{i_{k} + 1})^{e_{k-1} +e_{k}}  (x_{i_k} x_{i_{k} + 3})^{e_{k+1} -e_{k}} (x_{i_k +2} x_{i_{k} + 3})^{e_{k}} & \text{if } e_k \le e_{k+1}.
\end{cases}
\end{align*}
Read as a string, it satisfies Conditions (a) - (d). 

Since $(i_{k}+1, 2) = (i_{k+1}, j_{k+1}) < (i_{k+2}, j_{k+2})$ by assumption on Presentation \eqref{eq:initial rewriting}, we get 
$(i_{k} + 2, 1) \le (i_{k+2}, j_{k+2})$. It this is a strict inequality, then the string
\[
\bm' (x_{i_{k}+2} x_{i_{k+2} + j_{k+2}})^{e_{k+2}}
\]
satisfies Conditions (a) - (d). Otherwise, we have $(i_{k+2}, j_{k+2}) = (i_{k} + 2, 1)$, and we rewrite the monomial 
$\bm' (x_{i_{k}+2} x_{i_{k+2} + j_{k+2}})^{e_{k+2}}$ as the string 
\begin{align*}
& \cdots (x_{i_{k-2}} x_{i_{k-2} + j_{k-2}})^{e_{k-2}} \cdot \\
& \begin{cases}
(x_{i_k} x_{i_{k} + 1})^{e_{k-1} + e_{k+1}}  (x_{i_k} x_{i_{k} + 2})^{e_k -e_{k+1}} (x_{i_k +2} x_{i_{k} + 3})^{e_{k+1}+e_{k+2}}  & \text{if } e_k \ge e_{k+1},\\
(x_{i_k} x_{i_{k} + 1})^{e_{k-1} +e_{k}}  (x_{i_k} x_{i_{k} + 3})^{e_{k+1} -e_{k}} (x_{i_k +2} x_{i_{k} + 3})^{e_{k}+e_{k+2}} & \text{if } e_k \le e_{k+1}.
\end{cases}
\end{align*}
Read as a string, this description of $\bm' (x_{i_{k}+2} x_{i_{k+2} + j_{k+2}})^{e_{k+2}}$ satisfies Conditions (a) - (d). Repeating this rewriting step if necessary, eventually we can write $\bm$ in the required form. 

Conversely, it is clear that any string as described in the statement is the string representation of some monomial in $A$. 
\end{proof}

We use the above string presentation of a monomial to relate it to word in a suitable language.

\begin{defn}
    \label{def:language_counterexample}
(i) Define a language $\cL\subseteq \Sigma^*$ on the three-letter alphabet $\Sigma = \{\tau,\alpha_1,\alpha_2\}$  as the set of  words of the form 
\begin{align}
    \label{eq:language_counterexample}
    \tau^{k_1}\alpha_{i_1}\tau^{k_2}\alpha_{i_2}\dots \alpha_{i_d}\tau^{k_{d+1}}
\end{align}
with $d \in \N_0$ such that $k_1,\dots,k_{d+1}\geq 0$, $i_{1},\dots, i_{d}\in \{1,2\}$, and 
\begin{enumerate}[(a)]
    \item if $k_\nu=0$ for $1 < \nu \leq d$, then $i_{\nu-1} \le i_\nu$, 
    \item if $i_\nu=2$ and $k_{\nu+1}=1$ for some $\nu < d$, then $i_{\nu+1}=1$.
\end{enumerate}

(ii) 
Denote by $\cL_n^d$ the collection of words in $\cL$ with exactly $n$ occurrences of $\tau$ and $d$ occurrences of $\alpha_1$ and $\alpha_2$. 
\end{defn}

\begin{lm} 
    \label{lem:bijections}
For any $n\in \N$ and $d\in \N_0$, there is a bijection $[\Mon(A_n)]_d\rightarrow \cL_{n-1}^d$.
\end{lm}

\begin{proof}
Using the string presentation given in \Cref{prop:normal form}, this follows as in \Cref{lem:bijection}. We leave the details to the interested reader. 
\end{proof}

We now establish the main result of this section. 

\begin{thm} The equivariant Hilbert series of the filtration $\sA$ is 
\begin{align*}
    \equivH_\sA(s,t)&=\dfrac{ts^2+s}{-t^2s - ts^2 + t^2 + ts - 2t - s + 1}.
\end{align*}
\end{thm} 

\begin{proof} 
One checks that 
the language $\cL$ in \Cref{def:language_counterexample} is recognized by the following automaton: 
\begin{figure}[h]
    \centering
\begin{tikzpicture}[shorten >=1pt,node distance=5cm,auto]
  \tikzstyle{every state}=[fill={rgb:black,1;white,10}]
   \node[state,initial,accepting]           (1)     {$1$};
   \node[state,accepting]           (2) [below right of=1]     {$2$};
  \node[state,accepting] (3) [below left of=1]  {$3$};
  \path[->]
    (1)   edge  [loop right]         node {$\tau,\alpha_1$} (1)
     (2)   edge  [loop right]         node {$\alpha_2$} (2)
   (1)   edge  []         node {$\alpha_2$} (2)
    (2)   edge  []          node {$\tau$} (3)
    (3)   edge  [bend right]          node {$\alpha_1$} (1)
     (3)   edge  []          node {$\tau$} (1);      
\end{tikzpicture}
\caption{A finite automaton for the language in \Cref{def:language_counterexample}.}
\end{figure}
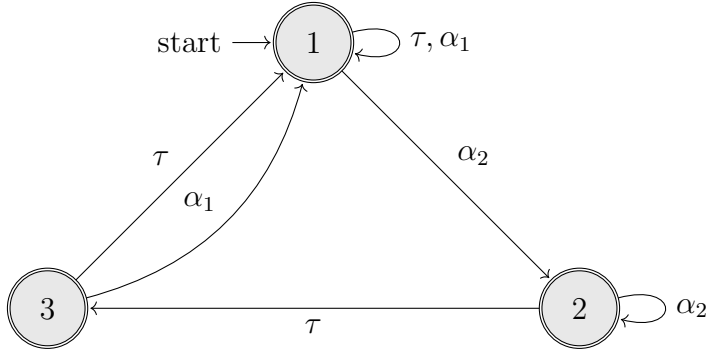

Thus, $\cL$ is a regular language. Using the weight function $\rho \colon \Sigma^* \to \K[s, t]$ with $\rho (\alpha_1) = \rho (\alpha_2) = t$ and $\rho (\tau) = s$ as well as \Cref{eq:automaton-HS}, one computes for the generating function,  
\begin{align*}
P_{\cL, \rho} (t,s)&=\begin{bmatrix}
    1\\1\\1
    \end{bmatrix}^T\begin{bmatrix}
    1-t-s &0 &-t-s \\
      -t &1-t &0\\
     0 &-s &1 
    \end{bmatrix}^{-1}\begin{bmatrix}
    1\\0\\0 
    \end{bmatrix}\\ \\
    &= \dfrac{ts+1}{-t^2s - ts^2 + t^2 + ts - 2t - s + 1}.
\end{align*}
Taking into account \Cref{lem:bijections}, it follows $\equivH_{\sA}(s,t)=s \cdot P_{\cL,\rho}(s,t)$. 
\end{proof}

\begin{cor}
   \label{cor:Hilb of non-fg ideal} 
The equivariant Hilbert series of $I$ is   
\begin{align*}
&\equivH_I (s,t)  =\sum\limits_{n\geq 1,d\geq 0} \dim_{\K} [I_n]_dt^ds^n \\
& = \dfrac{ts^4-ts^3-t^2s+s^3+ts-s}{t^2s^3+ts^4-t^3s-4t^2s^2-3ts^3+t^3+4t^2s+5ts^2+s^3-2t^2-6ts-3s^2+3t+3s-1}. 
\end{align*}
\end{cor} 

\begin{proof}
Since $A_n \cong \K[X_{[2] \times [n]}]/I_n$ we get 
\begin{align*}
\equivH_I (s,t) & =\sum\limits_{n\geq 1,d\geq 0} \dim_{\K} [I_n]_d t^ds^n  
= \sum\limits_{n\geq 1,d\geq 0} \big [\dim_{\K} [\K[X_{[2] \times [n] } ] ]_d  - \dim_{\K} [A_n]_d \big ] t^ds^n   \\
& =  \frac{(1-t)^2}{(1-t)^2 - s} - 1 - \equivH_{\sA} (s,t)
\end{align*} 
because 
\[
\sum\limits_{n\geq 1,d\geq 0} \dim_{\K} [\K[X_{[2] \times [n]}]]_d  t^ds^n  =  \frac{(1-t)^2}{(1-t)^2 - s} - 1
\]
(see, e.g., \cite[Example 6.1]{NR} or \cite[Proposition 2.6]{N}). A computation gives the claim. 
\end{proof}

We conclude with some comments on possible directions for further investigations. 
\smallskip 

\textbf{Future Directions.} This paper initiates the investigation of shift invariant monomial algebras and their presentation ideals 
and demonstrates that  they provide new interesting phenomena in representation stability. It is open to what extent our results generalize to arbitrary shift invariant monomial algebras or even shift invariant non-monomial algebras. 

This article utilizes connections between formal languages from computer science and infinite-dimensional objects in algebra. 
We introduce the Segre product of formal languages. It is worth investigating properties of this Segre product more systematically. The fact that the Segre product of regular languages is again a regular language could be of interest in other contexts as well. 

We used formal languages to establish rationality of equivariant Hilbert series of some filtrations. Hilbert's classical result includes a description of the denominator of Hilbert series of noetherian standard graded algebras (see, e.g.,\cite{HU}). For filtrations of algebras or modules, information on the rational functions appearing as equivariant Hilbert series have been obtained  in several papers, see, e.g., \cite{NR, N, SS-14}. Further work is needed. For example,  information on equivariant Hilbert series of filtrations determined by hierarchical models as investigated in \cite{MN} would be of interest. Note that we did not  analyze the structure of the used finite automata in this paper. Such an analysis could lead to new insights.  

%While we presented the example of a not finitely generated ideal up to the shift that has a rational equivariant Hilbert series, it is yet to be determined if finitely generated ideal implies rationality. 
%
%\begin{ques}
%For any $n,m\in \N$ take the polynomial ring $A_{n,m}=\K[x_{i,j} \mid i\in [n], j\in [n]]$ and $\sA=(A_{n,m})_{n,m\in \N}$, with finitely generated presenting ideal $I_{m,n}=\langle 0\rangle\subset A_{n,m}$. Is the equivariant Hilbert series 
%\[\equivH_{\sA}(t,s_1,s_2)=\sum_{d\geq 0, n,m\geq 1}\binom{mn+d-1}{d} t^d s_1^n s_2^m = \sum_{n,m\geq 1}\dfrac{1}{(1-t)^{nm}} s_1^n s_2^m  \]
%rational? 
%\end{ques}

\section*{Acknowledgments}
We thank  Corentin Bodart and Dietrich Kuske for sharing their comments and ideas on an earlier version of this paper. 
%Dietrich Kuske suggested the idea to prove  \Cref{thm:regularity of Segre language}. 

%\input{paper/Counterexample-aida}
%%%%%%%%%%%%%%%%%%%%%%%%%%%%%%%%

%\bibliographystyle{plain}
%\bibliography{Paper/references}

%\hfill
%
%\noindent
%\footnotesize {\bf Authors' addresses:}
%
%\hfill
%
%\noindent Aida Maraj\\
%Department of Mathematics, \\ Universiry of Michigan, Ann Arbor (MI), United States of America\\
%\hfill {\tt maraja@umich.edu}
%
%\hfill
%
%\noindent Uwe Nagel\\
%Department of Mathematics, \\University of Kentucky, Lexington (KY), United States of America \\
%\hfill {\tt uwe.nagel@uky.edu}

\end{document}